\documentclass[a4paper,10pt]{article}

\author{Léo Brunswic}
\def\title{Cauchy-Compact flat spacetimes with BTZ singularities}

\usepackage[utf8]{inputenc}
  \usepackage[usenames,dvipsnames]{xcolor}
  \usepackage{amsthm}
  \usepackage{yhmath,amsfonts,amssymb,amsmath,MnSymbol, mathrsfs,amsbsy}
  \usepackage{mdframed}
  \usepackage[all]{xy}
  \usepackage[english]{babel}
  \usepackage[colorlinks=true, linkcolor=black,citecolor=black]{hyperref}
  \usepackage{fullpage, version}
  \usepackage{geometry}
  \usepackage[toc, page]{appendix} 
  \usepackage{graphicx}
  \usepackage{enumerate}
  \usepackage{tkz-euclide}
  \usepackage{ulem}
  \usepackage{color}

  \usetkzobj{all}
\def\MM{\overline M}
  \def\Z{\mathbb{Z}}
  \def\H{\mathbb{H}}
  \def\M{\mathcal{M}}

  \def\W{\mathcal W}
  
  \def\stab{\mathrm{Stab}}

  \def\U{\mathcal U}
  
  \def\Hom{\mathrm{Hom}}
  \def\teich{\mathrm{Teich}}
\def\V{\mathcal V}

\def\C{\mathscr C}
\def\R{\mathbb{R}}
\def\E{\mathbf E}

\def\N{\mathbb N}
\def\T{\mathcal{T}}

\def\d{\mathrm{d}}

\def\D{\mathcal{D}}
\def\SO{\mathrm{SO}}
\def\E{\mathbb{E}}
\def\BTZ{{\E_{0}^{1,2}}}
\def\sing{\mathrm{Sing}}
\def\susp{\mathrm{susp}}

\def\isom{\mathrm{Isom}}

\def\BTZext{\mathrm{BTZ-ext}}
\def\reg{\mathrm{Reg}}

\def\fix{\mathrm{Fix}}
\def\dec{\mathrm{dec}}

\newcommand{\fonction}[5]{\displaystyle#1:\begin{array}{l|rcl}
& \displaystyle #2 & \longrightarrow & \displaystyle #3 \\
    & \displaystyle #4 & \longmapsto & \displaystyle #5 \end{array}}
\newcommand{\fonctionn}[4]{\displaystyle\begin{array}{l|rcl}
& \displaystyle #1 & \longrightarrow & \displaystyle #2 \\
    & \displaystyle #3 & \longmapsto & \displaystyle #4 \end{array}}

\newtheorem{theo}{\bf{Theorem}}

\newtheorem{theo_ext}{\bf{Theorem}}[section]

\newtheorem{lem}{Lemma}[section]
\newtheorem{cor}[lem]{Corollary}
\newtheorem{prop}[lem]{Proposition}
\newtheorem{defi}[lem]{Definition}

\newtheorem{rem}[lem]{Remark}

\usepackage{float}
\begin{document}
 
\thispagestyle{empty}
\begin{minipage}[t]{0.5\textwidth}
\textit{Laboratoire de Mathématiques d'Avignon \\
Université d'Avignon et des pays du Vaucluse}
\end{minipage}
\hfill
\begin{minipage}[t]{0.35\textwidth}
\textbf{Léo Brunswic\\
\today\\
}
\end{minipage}
\vspace{.5cm}
\begin{center}
\begin{minipage}{0.7\textwidth}
\hrule \vspace*{1.0cm}
\begin{center}
\LARGE
\textbf{\title }
\vspace{.5cm}
\end{center}
\hrule
\end{minipage}
\vspace{4cm}
\end{center}

\begin{center}
\hspace{1.5cm}

\end{center}

\vspace{-0cm}
\begin{center}

\end{center}

\def\Hol{\mathrm{Hol}}
\def\so{\mathfrak{so}}
\def\pen{\mathscr{P}}

\begin{abstract}
 The zoology of singularities for Lorentzian manifold is slightly more complex than for Riemannian manifolds. 
 Our present work  
 study Cauchy-compact globally hyperbolic singular flat spacetimes with extreme BTZ-like singular lines. We use 
  the notion of BTZ-extension of a singular spacetime introduced in a previous paper to give a description of 
  Moduli spaces of such manifolds in term of common Teichmüller spaces.   
  This description is used to construct convex polyhedral cauchy-surface in Cauchy-compact flat spacetimes with BTZ.

\end{abstract}
\newpage
	\setcounter{tocdepth}{3}
	\tableofcontents
	 
\section{Introduction}
   \subsection{Context and motivation}
      
      Let $\E^{1,2}$ be Minkoswki space, namely $\R^3$ together with the quadratic form $\d s^2= -\d t^2+\d x^2+\d y^2$, 
      and let  $\isom(\E^{1,2})=\SO_0(1,2)\rtimes \R^3$ be its direct time-orientation preserving isometry group. 
      The main object of this paper are {\it Cauchy-complete Cauchy-maximal globally hyperbolic singular } $(\isom(\E^{1,2}),\E^{1,2} )$-manifolds.      
      The most simple example of such a $\E^{1,2}$-manifold $M$ is given by a quotient of $F= \{-t^2+x^2+y^2<0, t>0\}$
      by a Fucshian group, say $\Gamma= \Gamma(2)$ the index  6 congruence subgroup of $\mathrm{SL}(2,\Z)$. 
      The subset $\{-t^2+x^2+y^2 = -1, t>0\}$ is a natural equivariant embedding of $\H^2$ into $F$ giving a natural 
      {\it Cauchy-surface} of $F$ and thus of $M:=F/\Gamma$.
      Consider a $\Gamma$-invariant triangulation of $\H^2$, say 
      $\{\gamma T_i : \gamma \in \Gamma, i \in \{1,\cdots,6\}\}$, 
      we can take the suspension of each $T_i$,  $\susp(T_i) := (\R_+^* \times T_i, -\d t^2+t^2\d s_{T_i}^2 )$  
      and glue these cones face to face to re-construct $M$. In this constuction, the lightlike edges of the suspension 
      have been implicitely removed. If we extend the gluing to the lightlike edges 
      of the suspension, we obtain a $\E^{1,2}$-manifold with {\it extreme BTZ-like singular lines} that is 
      a $\E^{1,2}_0$-manifold.
      We thus constructed a BTZ-extension \cite{brunswic_btz_ext} of $M$ say $M'$.
      It is easy to construct a polyhedral compact Cauchy-surface of $M'$ however contrary to the natural embedding of
      $\H^2/\Gamma$ in $M$, such polyhedral surface may not be convex. 
      A good construction of convex polyhedral Cauchy-surfaces 
      is given by the Penner-Epstein surface Penner  describes in \cite{penner1987,MR3052157}. 
      The construction can be rephrased using the singular spacetime terminology:
      the boundary of the closed convex hull of a family of one point chosen  on  each BTZ-line 
      is a convex polyhedral Cauchy-surface of $M'$.
      
      This paper generalises these constructions to Cauchy-maximal Cauchy-compact globally hyperbolic 
      $\E^{1,2}$-manifolds
      with BTZ or $\E^{1,2}_0$-manifolds following \cite{brunswic_btz_ext} terminology. The main result of \cite{brunswic_btz_ext} states that the regular part of such a manifold 
      is a Cauchy-complete Cauchy-maximal globally hyperbolic $\E^{1,2}$-manifold. However, 
      we don't {\it a priori} have a constant curvature
      Cauchy-surface in a general Cauchy-complete Cauchy-maximal $\E^{1,2}$-manifold. Therefore, a reconstruction from a 
      triangulation of a hyperbolic Cauchy-surface is difficult to obtain and we shall proceed in a different manner.
      More precisely, we avoid this difficulty by 
      constructing the maximal BTZ-extension of Cauchy-complete Cauchy-maximal $\E^{1,2}$-manifold starting 
      from the description of its universal cover as a {\it regular  domain} in Minkoswki 
      given by Mess, Bonsante, Benedetti and Barbot \cite{mess, bonsante_flat_2003, barbot_globally_2004,MR2499272}.
    
         \subsection{Terminology and dependancies}
	This paper follows directly \cite{brunswic_btz_ext}, many elementary properties of  $(G,X)$-structures and 
	$\E^{1,2}_A$-manifolds that are useful to the present work can be found in this previous paper.
	
	We use freely basics elements of the theory of $(G,X)$-structures in particular holonomy and developping map. 
	See for instance section 4 of \cite{goldman_course_note}.   Most of the structures we use are special cases 
	of $\E^{1,2}_A$-manifolds 
	in the sense of \cite{brunswic_btz_ext}, i.e singular flat spacetimes. These are not $(G,X)$-manifolds, however they contain a regular locus i.e.
	a dense open subset which is a $(\isom(\E^{1,2}),\E^{1,2})$-manifold. The complement of the regular locus is called the singular locus. 
	The developping map and holonomy of a $\E^{1,2}_A$-manifold are defined as the  developping map and the holonomy of
	its regular locus.
	By definition, points of the singular locus are locally modeled on a singular model space. In most of the paper,
	$A=\{0\}$ : the local model spaces of the singular spacetimes we are handling are $\E^{1,2}$ and $\E^{1,2}_0$.
	Details on the geometry of the BTZ model space $\E^{1,2}_0$ can be found in section 1 of \cite{brunswic_btz_ext} .
	$\E^{1,2}_A$-manifolds are singular spacetimes, the theory of (regular) semi-riemannian spacetimes is detailed in  \cite{oneil}.
	Such spacetimes come with two natural orders called the time order and the causal order.
	Causal (resp. timelike) curves are  contiuous monotonic curves for the causal (resp. time) order.
	For a detailed exposition of properties of causal orders in semi-riemannian spacetimes see \cite{MR2436235}. 
	Many fundamental definitions and results extend to singular spacetimes \cite{Particles_1,Particles_2, brunswic_btz_ext}. 
	We will freely use fundamental definition and results about globally hyperbolic manifolds such as 
	diamonds, Cauchy-surface, Cauchy-maximality, Cauchy-completeness and time functions \cite{oneil,MR0270697,MR0250640,andersson_time_function}.

	A Cauchy-complete flat spacetime is absolutely maximal if it cannot be strictly embedded in any other 
	Cauchy-complete flat spacetime. Properties of such spacetimes are given in Section 2.2 and extended to Cauchy-complete flat spacetime with BTZ in 
	section 3.1.4.
	
	Properties about Teichmüller theory will be given when needed. Most of the results we use can be found in 
	\cite{MR3052157}.
    
      \subsection{Results}
      The main results of the paper are :
      \setcounter{theo}{1}
      \begin{theo}
       	Let $\Sigma$ be a closed surface and $S$ a finite subset. 
       There is a canonical identification between the tangent fiber bundle of the Teichmüller space of $\Sigma\setminus S$ and 
       the moduli space of marked Cauchy-maximal Cauchy-compact globally hyperbolic spacetimes on $\Sigma\times \R$
      with $|S|$  BTZ-lines
       \end{theo}
       \setcounter{theo}{3}
      \begin{theo} Given a Cauchy-maximal Cauchy-compact globally hyperbolic $\E^{1,2}_0$-manifolds and 
      an arbitrary choice of points $p_1,\cdots, p_s$  in every BTZ-line $\Delta_1,\cdots \Delta_s$, 
      there exists a unique convex polyhedral Cauchy-surface of vertices $p_1,\dots,p_s$.
      \end{theo}
      
      Secondary results are worth noting : 
      \setcounter{theo}{0}
      \begin{theo} A description of the maximal BTZ-extension of a Cauchy-maximal Cauchy-complete globally hyperbolic $\E^{1,2}$-manifold
       as the quotient of a convex domain of Minkoswki space. 
      \end{theo}
         \setcounter{theo}{2}
      \begin{theo} The moduli space of marked singular Euclidean structures on a marked surface $(\Sigma,S)$
       is canonically identified with the decorated Teichmüller space of 
       $(\Sigma,S)$ and the marked moduli space of linear Cauchy-compact $\E^{1,2}_0$-structures on $(\Sigma,S)$.
       \end{theo}
          \setcounter{theo}{0}

   \section{Flat Lorentzian manifolds and Moduli spaces}\label{sec:special_introduction}
\subsection{ Teichmüller space}

Let $\Sigma$ be a compact surface of genus $g$ and let $S$ be a finite subset of cardinal $s\geq 0$ such that $2g-2+s >0$. 
The punctured surface 
$\Sigma^* :=\Sigma\setminus S$ can then be endowed with a complete hyperbolic metric of finite volume.
A marked surface homeomorphic to $\Sigma^*$ is a couple $(\Sigma_1^*,h_1)$ where $\Sigma_1^*$ is a surface  and $h_1 
:\Sigma^*  \rightarrow \Sigma_1^*$ is a homeomorphism.
Two marked complete hyperbolic surfaces of finite volume $(\Sigma_1^*,h_1,m_1)$ and $(\Sigma_2^*,h_2,m_2)$, where $m_i$ is a hyperbolic metric on 
$\Sigma_i^*$,
homeomorphic to $\Sigma^*$ are equivalent if there exists an isometry $\varphi : (\Sigma_1^*,m_1)\rightarrow (\Sigma_2^*,m_2)$ 
such that $h_2^{-1}\circ\varphi \circ h_1$ is a homeomorphism of $\Sigma^*$ homotopic to the identity on $\Sigma^*$.

$$(\Sigma_1^*,m_1,h_1) \sim (\Sigma_2^*,m_2,h_2) \quad \Leftrightarrow \quad \begin{minipage}{5cm}
  \xymatrix{&(\Sigma_1^*,m_1) \ar@{-->}[dd]^{\exists\varphi}\\ \Sigma^* \ar[ur]^{h_1} \ar[dr]_{h_2}&\\& (\Sigma_2^*,m_2) }
  \end{minipage} ~\mathrm{with}~\left\{\begin{array}{l}
h_2^{-1}\circ \varphi \circ h_1 \sim \mathrm{Id}_{\Sigma^*}  \\  \varphi ~\mathrm{isometry}\end{array}\right.$$
The Teichmüller space of $(\Sigma,S)$, denoted $\teich_{g,s}$, 
is the space of all complete hyperbolic marked surface  of finite volume homeomorphic to $\Sigma^*$ up to equivalence.  
Let $\Gamma:=\pi_1(\Sigma^*)$,
a point of Teichmüller space $\teich_{g,s}$, as all $(G,X$)-manifolds \cite{ratcliff_foundation,goldman_course_note}, has an holonomy which is a point  
of $\Hom(\Gamma,\SO_0(1,2))/\SO_0(1,2)$ where $\SO_0(1,2)$ acts by conjugacy.  
This defines an injective map \cite{ratcliff_foundation,MR3052157}  $\Hol :\teich_{g,s} \rightarrow  \Hom(\Gamma,\SO_0(1,2))/\SO_0(1,2)$.
The image is the so-called Teichmüller component \cite{MR1174252}  of $\Hom(\Gamma,\SO_0(1,2))/\SO_0(1,2)$ and can be described in the following way.
$\Gamma$ has the following presentation
$$\Gamma = \left\langle a_1,b_1,\cdots,a_g,b_g, c_1,\cdots,c_s \left |  \prod_{i=1}^g [a_i,b_i] \prod_{j=1}^s c_j=1\right.\right\rangle.$$
The generators $c_i$ are called peripherals and correspond to loops around the punctures $S$.  
The holonomy $\rho$ of a point of $\teich_{g,s}$ is marked by a choice of such generators
$(a_i,b_i,c_j)_{i\in [1,g] ; j\in[1,s]}$ of $\Gamma$. 
A marked linear representation $\rho : \Gamma\rightarrow \SO_0(1,2)$ is in the image of  $\teich_{g,s}$ by $\Hol$
if and only if it is admissible in the following sense.
\begin{defi}[Admissible representation (linear case)]\label{defi:admissible_linear}\ \\
Let $\Gamma = \left\langle a_1,b_1,\cdots,a_g,b_g, c_1,\cdots,c_s \left |  \prod_{i=1}^g [a_i,b_i] \prod_{j=1}^s c_j=1\right.\right\rangle$
be a marked surface group.
A marked representation $\rho:\Gamma \rightarrow \SO_0(1,2)$ is admissible if 
\begin{itemize}
 \item $\rho$ is faithful and discrete;
 \item for all $j\in \{1,\cdots, s\}$, $\rho(c_j)$ is parabolic;
 \item for all $i\in\{1,\cdots, g\}$, $\rho(a_i)$ and $\rho(b_i)$ are hyperbolic.
\end{itemize}
 \end{defi}
\begin{rem} This is exactly the definition of a lattice of $\SO_0(1,2)$.
 
\end{rem}
 $\Hom(\Gamma,\SO_0(1,2))/\SO_0(1,2)$ has the natural compact-open topology and the Teichmüller component is a differential manifold 
 diffeomorphic to $\R^{6g-6+2s}$. In the following we identify $\teich_{g,s}$ and the Teichmüller component of  $\Hom(\Gamma,\SO_0(1,2))/\SO_0(1,2)$.  
 
  We now give a  description of the tangent fiber bundle of $\teich_{g,s}$ following Goldman \cite{MR762512}.
  Let $[\rho] \in \teich_{g,s}$ be a class of marked representation. The tangent space to $\Hom(\Gamma,\SO_0(1,2))$ at $\rho$ is 
  the space of cocycles for $\rho$, i.e. the space of $\tau : \Gamma \rightarrow \mathfrak{so}(1,2)$ such that
   $$\forall \gamma_1,\gamma_2, \tau(\gamma_1\gamma_2)=\tau(\gamma_1)+\mathrm{Ad}(\rho(\gamma_1))\tau(\gamma_2). $$
  Moreover, the action of $\SO_0(1,2)$ by conjugacy induces an equivalence of cocycles via coboundary i.e. the cocycles $\tau:\Gamma \rightarrow \so(1,2)$
  for which there exists $u\in \so(1,2)$ such that 
  $$\forall \gamma \in \Gamma, \tau(\gamma)=\mathrm{Ad}(\rho(\gamma))u-u.$$
Then, for $\rho : \Gamma \rightarrow \SO_0(1,2)$ admissible, 
the tangent space $T_{\H^2/\rho}\teich_{g,s}$ is naturally identified with 
a subspace of $H^1_{\mathrm{Ad}\circ\rho}(\Gamma, \so(1,2))$. 
 
Consider $j\in [1,s]$ and a 1-parameter family $(\rho_s)_{s\in\R}$ of admissible representations with $\rho_0=\rho$. The image $\rho_s(c_j)$ is parabolic 
for all $s\in \R$ thus there exists a 1-parameter family $(\phi_s)$ of elements of $\SO_0(1,2)$ such that for all $s\in\R$, 
$\rho_s(c_j)=\phi_s\rho(c_j)\phi_s^{-1}$. A simple computation shows there exists $u$ such that
$$\left.\frac{\d \rho_s}{\d s}\right |_{s=0}(c_j)=\mathrm{Ad}(\rho(c_j))u-u$$
 Thus, let $\tau$ be a tangent vector at $\rho$, for $j\in[1,s]$, $\tau(c_j)$ is orthogonal to the line of fixed points of 
$ \mathrm{Ad}(\rho(c_j))$;
 the orthogonal being taken for the Killing form on $\so(1,2)$. When $s>0$, $\Gamma$ is a free group and thus dimensions are 
easily computed and show that the tangent vectors at $\rho$
are exactly the cocycles satisfying this property up to coboundaries. This is still true when $s=0$, see \cite{MR762512}.

Notice that the Killing form on $\so(1,2)$ is non-degenerate of signature $(1,2)$ thus $\so(1,2)$ is isometric to $\E^{1,2}$.
Furthermore, adjoint action of $\phi\in\SO_0(1,2)$ on $\so(1,2)$ is hyperbolic (resp. parabolic, resp. elliptic) if and only if $\phi$ is hyperbolic
(resp. parabolic, resp. elliptic). A point of $T\teich_{g,s}$  can then be seen as a marked representation $\Gamma \rightarrow \isom(\E^{1,2})$ up to conjugacy.
\begin{defi}\label{defi:admissible} Write $L$ the projection $\isom(\E^{1,2})\rightarrow \SO_0(1,2)$. Let  $\rho:\Gamma\rightarrow \isom(\E^{1,2})$ be 
a representation.
\begin{itemize}
 \item The linear part of $\rho$, is $\rho_L:=L\circ \rho$.
 \item The cocycle part of a $\rho$, is $\tau_\rho:=\rho-\rho_L$.
\end{itemize}

For $\phi\in \isom(\E^{1,2})$, we also write $\phi_L=\phi_L$ and $\tau_\phi := \phi-\phi_L$. 
\end{defi}
\begin{defi} For $\phi\in \isom(\E^{1,2})$, $\fix(\phi)  = \{p\in \E^{1,2} ~|~ \phi x=x\}$ is a fixator of $\phi$.
 
\end{defi}

\begin{defi} Let $\phi \in \isom(\E^{1,2})$, $\tau_\phi$ is tangent if $\tau_\phi \in \fix\left(\rho(c_j)\right)^{\perp}$.
 
\end{defi}

The tangent bundle of $\teich_{g,s}$ is then the set of admissible representations in the following sense.
\begin{defi}[Admissible representation (affine case)] \ \\ Let $\Gamma = \left\langle a_1,b_1,\cdots,a_g,b_g, c_1,\cdots,c_s \left |  \prod_{i=1}^g [a_i,b_i] \prod_{j=1}^s c_j=1\right.\right\rangle$
be a marked surface group.
A marked representation $\rho:\Gamma \rightarrow \isom(\E^{1,2})$ is admissible if 
\begin{itemize}
 \item $\rho_L$ is admissible;
 \item $\tau_\rho(c_j)$ is tangent for every $j$.
\end{itemize}
\end{defi}

\begin{prop}\label{prop:teich_tangent}
The tangent fiber bundle of Teichmüller space  $T \teich_{g,s}$  is 
canonically identified with the set conjugacy classes of marked representations $\Gamma\rightarrow \isom(\E^{1,2})$.
\end{prop}

\subsection{Globally hyperbolic Cauchy-complete flat spacetimes}

We briefly give main results about regular domains, a more complete study is given in 
\cite{bonsante_flat_2003}, \cite{barbot_globally_2004} and \cite{MR2499272}. 
We use notations from \cite{oneil},  the timelike (resp. causal) future of a set $A$ is denoted $I^+(A)$ (resp. $J^+(A)$).
Such a manifold has a natural order relation, the causal order.  The map 
$$\fonction{rev}{\E^{1,2}}{\E^{1,2}}{(t,x,y)}{(t,x,y)}$$ 
is the time reversal map. This transformation preserves the quadratic form, it is in $\mathrm{O}(1,2)$ but not in $\SO_0(1,2)$. 
The time reversal map induces an involution on the set of $\E^{1,2}$-manifolds : let $M$ be a $\E^{1,2}$-manifold, 
we can replace every local chart $\U\subset M$, $\varphi:\U\rightarrow \E^{1,2}$ by $rev\circ \varphi : \U \rightarrow \E^{1,2})$.
This inverse the causal order on $M$. This transformation is called {\it time reversal}.

An affine lightlike plane  in $\E^{1,2}$ 
is of the form $\{x\in \E^{1,2},\langle x|u\rangle=\lambda\}$ for some $u$ lightlike vector and $\lambda\in \R$. The set of affine 
lightlike plane is then  homeomorphic to $S^1\times \R$.

	\begin{defi}[Regular domain] \label{def:regular_domain}
	  A regular domain $\Omega\subset \E^{1,2}$ is a set of the form 
	  $$\Omega_\Lambda:=\bigcap_{\Pi \in \Lambda} I^+(\Pi)$$
	  for some closed family of lightlike planes $\Lambda$.
	\end{defi}

	\begin{theo_ext}[\cite{barbot_globally_2004}]\label{barbot} Let $\Gamma$ be a discrete subgroup of $\isom(\E^{1,2})$ torsionfree such that 
	$\Gamma \cdot \Omega =\Omega$. Then $\Gamma$ acts properly discontinuously on $\Omega$ and 
	$\Omega/\Gamma $  is a globally hyperbolic Cauchy-complete $\E^{1,2}$-manifold.

	\end{theo_ext}
      
	\begin{theo_ext}[\cite{barbot_globally_2004,MR2499272}]\label{theo:barbot_mess}
	  Let $M$ be  a maximal globally Cauchy-complete hyperbolic spacetime. Then the developping map is an 
	  embedding so that the universal covering $\widetilde M$ of $M$ can be identified with a domain 
	  $\Omega\subset \E^{1,2}$. Moreover, up to time inversion, one of the following holds.
	  \begin{enumerate}[(i)]
	  \item $\Omega = \E^{1,2}$ and the holonomy group acts as a free abelian group of spacelike translations of rank at most 2.
	  \item $\Omega$ is the future of a lightlike plane $\Pi$ and the holonomy group, if non
	  trivial, is generated by a spacelike translation or a parabolic linear isometry. 
	  \item $\Omega = I^+(\Pi^-)\cap I^-(\Pi^+)$, where $\Pi^+,\Pi^-$ are parallel lightlike planes. The holonomy group, if non
	  trivial, is generated by a spacelike translation or a parabolic linear isometry. 
	  \item $\Omega$ is a regular domain, the linear part of the
	  holonomy $\pi_1(M) \rightarrow \SO_0(1,2)$ is a faithfull and discrete representation.
	  \end{enumerate}
      \end{theo_ext}
      \begin{cor} In particular a Cauchy-maximal globally hyperbolic Cauchy-complete spacetime is future complete up to time reversal.       
      \end{cor}

      \begin{defi}
      A spacetime  $M$ is of type (i) (resp. (ii), resp. (iii), resp (iv)) if  it is a globally hyperbolic,
      Cauchy-complete $\E^{1,2}$-manifold falling into case (i) (resp. (ii), resp. (iii), resp. (iv))  of Theorem \ref{theo:barbot_mess}
      \end{defi}

\begin{defi}
    Let $M$ be  a globally hyperbolic Cauchy-complete $\E^{1,2}$-manifold, $\Gamma$ be its fundamental group and $\rho$ its holonomy. 
 $M$ is absolutely maximal if for all globally hyperbolic Cauchy-complete $\E^{1,2}$-manifold $M'$ and
 for all injective morphism of $\E^{1,2}$-structure $i:M\rightarrow M'$, $i$ is surjective.
\end{defi}

\begin{prop}[\cite{barbot_globally_2004}] \label{prop:abs_ext}
 Let $M$ be a globally hyperbolic Cauchy-complete $\E^{1,2}$-manifold, then there exists an 
absolutely maximal globally hyperbolic Cauchy-complete $\E^{1,2}$-manifold $\overline M$ in which $M$ embeds isometrically. Moreover, $\overline M$ is
unique up to isomorphism.
\end{prop}

\begin{prop}[\cite{barbot_globally_2004}] \label{prop:abs_ext_2}
Let $M$ be a globally hyperbolic Cauchy-complete $\E^{1,2}$-manifold, 
let $\Gamma:=\pi_1(\M)$ and let $\rho: \Gamma\rightarrow \isom(\E^{1,2})$ be its holonomy. 
  Define $$\Omega(\rho) = \bigcap _{\Pi \in \Lambda(\rho)} I^+(\Pi)$$ where  $\Lambda(\rho)$ is the closure of the 
  set of lightlike plane 
  which are repellent fixed point of hyperbolic element of $\rho(\Gamma)$. 
  
 Then the absolutely maximal extension of $M$ is isomorphic to $\Omega(\rho)/\rho$.
\end{prop}

\begin{lem} Let $\phi \in \isom(\E^{1,2})$ parabolic, the following are equivalent :
\begin{enumerate}
 \item[(i)]  $\mathrm{Fix}(\phi)\neq \emptyset$
 \item[(ii)] $\tau_\phi \in \mathrm{Fix}(\phi_L)^\perp$
\end{enumerate}
\end{lem}
\begin{proof}
 To begin with, these two properties are invariant under conjugation of $\phi$.
  Assume $\fix(\phi)\neq \emptyset$ and let $p\in \E^{1,2}$ such that $\phi p = p$, up to a conjugation 
  we can assume $p=O$ and thus $\tau_\phi=0 \in \fix(\phi_L)^\perp$.
  
  Assume $\tau_\phi\in \fix(\phi_L)^\perp$. Conjugating $\phi$ by a translation of vector $u$ changes
  $\tau_\phi$ into $\tau_\phi + (\phi_L-1)u$. The map $\E^{1,2} \xrightarrow{\phi_L-1} \E^{1,2}$ is 
  linear of rank 2 since $\phi_L$ has exactly one direction of fixed points. 
  Since $\mathrm{Im}(\phi_L-1) \subset \fix(\phi_L)^\perp$ and both are of dimension 2, then 
  $\mathrm{Im}(\phi_L-1) = \fix(\phi_L)^\perp$. Therefore, there exists $u \in\E^{1,2}$ such that 
  the conjugation of $\phi$ by the translation of vector $u$ is linear. Finally, $\phi$ is conjugated to 
  an isometry which admits a fixed point and thus $\fix(\phi)\neq \emptyset$. 

\end{proof}

\begin{lem}\label{lem:tangent}
 Let $\Omega$ be a regular domain stabilized by some torsionfree discrete subgroup $G\subset \isom(\E^{1,2})$. Then 
 for all $\phi\in G$ parabolic, $\tau_\phi$ is tangent. 
\end{lem}
\begin{proof}
 See \cite{barbot_globally_2004} section 7.3.
\end{proof}

\begin{cor} \label{cor:admissible} Let $\Sigma$ be a compact surface, $S$ be a finite subset of $\Sigma$ and $\Gamma=\pi_1(\Sigma^*)$ such that 
$2g-2+s>0$. 
Let $M$ be a globally hyperbolic Cauchy-complete spacetime homeomorphic to $(\Sigma^*) \times \R$
 and $\rho:\Gamma \rightarrow \isom(\E^{1,2})$ be its holonomy. 
 
 Then $\rho$ is admissible if and only if $\rho_L$ is admissible.
\end{cor}
\begin{proof}If $\rho$ is admissible, by definition, $\rho_L$ is admissible. If $\rho_L$ is admissible then $\rho(\Gamma)$
fixes some non-empty regular domain and thus, by Lemma \ref{lem:tangent}, $\tau_{\rho(\gamma)}$  is tangent whenever $\rho(\gamma)$
is parabolic and thus $\rho$ is admissible. 

\end{proof}

We end this section by a remark on the topology of Cauchy-complete globally hyperbolic flat spacetimes. 
\begin{prop} Let $\Sigma^*$ be a punctured surface of genus $g\geq 0$ with $s\geq 0$ punctures such that $2g-2+s>0$
and $\Gamma$ its fundamental group. Let $\rho : \Gamma \rightarrow \isom(\E^{1,2})$ be a discrete faithful 
marked representation of $\Gamma$.  Let $M$ be a globally hyperbolic Cauchy-complete spacetime of fundamental group $\pi_1(M)$ isomorphic to $\Gamma$.
 
 If the holonomy of $M$ is $\rho$ then $M$ is homeomorphic to $\Sigma^* \times \R$.
\end{prop}
\begin{proof}[Sketch of proof]\label{prop:topology}
We let $\Gamma$ acts on $\E^{1,2}$ via $\rho$ and acts on $\H^2$ via $\rho_L$. Since $\Gamma$ in not abelian, 
from Theorem \ref{theo:barbot_mess}, $\rho_L$ is faithful and discrete.

Write $\Omega$ a $\Gamma$-invariant regular domain such that $M$ is isomorphic to $\Omega/\Gamma$.
From Proposition 3.3.3 in \cite{MR2499272}, there exists a $\Gamma$-invariant convex $\C^1$ Cauchy-surface of $\Omega$, say 
$\widetilde \Sigma_1^*$ and $\Sigma^*_1 :=\widetilde \Sigma_1^*/\Gamma$ is a Cauchy-surface of $M$. 
The Gauss map of $\widetilde \Sigma_1^*$ defines a map $N: \Sigma^*_1   \rightarrow \H^2$.  
$\widetilde \Sigma_1^*$ is the graph of a convex function defined on $\R^2$ thus from a Theorem of Minty 
\cite{MR0132379,MR0167859}, there exists a closed convex domain $H\subset \H^2$ such that $Int(H) \subset N(\Sigma_1^*) \subset H$.
We now prove that $N$ is  proper. Sections 3.4 to 3.6 of \cite{MR2499272} explain how to construct a geodesic lamination on $H$ 
from $\Omega$ and Proposition 3.6.1 shows that the inverse image by $N$ of a compact curve $c$ in $H$ is of finite length if and 
only if it does not intersects $\partial H$. This implies that $N$ is proper on $N^{-1}(Int(H))$ ans thus 
$N^{-1}(Int(H))/\Gamma$ and $Int(H)/\Gamma$ have the same number of ends.  Since 
$Int(H)/\Gamma$ and $N^{-1}(Int(H))/\Gamma$ have the same finitely generated fundamental group and the same number of ends, they are homeomorphic.

A simple analysis of the way $\widetilde\Sigma_1^*$ is constructed shows that 
the inverse image of the boundary of $H$ is a set of the form $\gamma +a+\R_+\cdot h$ for some spacelike vector $h$, some vector $a$ and some 
 $\gamma$ geodesics of $\H^2$. Therefore, $\Sigma^*_1$ is obtained from $N^{-1}(Int(H))/\Gamma$ by extending its ends, 
 it is thus homeomorphic. Finally, $\Sigma^*_1$ is homeomorphic to $Int(H)/\Gamma$, which is homeomorphic to $\Sigma^*$.

\end{proof}

 \subsection{Flat Lorentzian Moduli spaces}

The definition of $\teich_{g,s}$ can be adapted to define Moduli spaces associated to marked $\E^{1,2}_A$-structures.
\begin{defi}[Equivalence of marked $\E^{1,2}_A$-manifolds]Let $\Sigma$ be a non necessarily compact surface.

 Let $(M_1,h_1,m_1)$, $(M_2,h_2,m_2)$ be marked $\E^{1,2}_A$-manifolds, 
 where $h_i: \Sigma\times \R \rightarrow M_i$ is a homeomorphism
 and $m_i$ is a  $\E^{1,2}_A$-structure on $M_i$.
 The manifolds $M_1$ and $M_2$  are equivalent if there exists a   $\E^{1,2}_A$-isomorphism  $\varphi : M_1 \rightarrow M_2$,
 such that $h_2^{-1}\circ \varphi \circ h_1$ is homotopic to $(x,t)\mapsto (x,t)$ or $(x,t)\mapsto (x,-t)$.
 
\end{defi}
 \begin{rem} As stated before, the reversal of time is a transformation among spacetimes which is not an isomorphism of $(\isom(\E^{1,2},\E^{1,2})$-structure.
 However, the properties of a $\E^{1,2}$-manifold is closely related to the properties of its time reversal.
 We thus authorise time reversal as equivalence.
 
\end{rem}

\begin{defi}[Linear marked $\E^{1,2}$-Moduli space] Let $\Sigma$ be a compact surface of genus $g$, $S$ be a finite subset of cardinal 
$s>0$.
  The Moduli space $\M_{g,s}^L(\E^{1,2})$ is the space of  equivalence classes of marked
 $\E^{1,2}$-manifolds of \textbf{linear admissible} holonomy homeomorphic to $(\Sigma^*)\times \R$
 which are  globally hyperbolic Cauchy-complete, Cauchy-maximal. 
\end{defi}

\begin{defi}[Marked $\E^{1,2}$-Moduli space]Let $\Sigma$ be a compact surface of genus $g$, $S$ be a finite subset of cardinal 
$s>0$.
The Moduli space $\M_{g,s}(\E^{1,2})$ is the space of equivalence classes of
 marked $\E^{1,2}$-manifolds of \textbf{affine admissible} holonomy homeomorphic to $(\Sigma^*)\times \R$
 which are  globally hyperbolic Cauchy-complete, Cauchy-maximal and \textbf{absolutely maximal}. 
 \end{defi}
 
\begin{defi}[Linear marked $\E^{1,2}_0$-Moduli space]Let $\Sigma$ be a compact surface of genus $g$, $S$ be a finite subset of cardinal 
$s>0$.
 The Moduli space $\M_{g,s}^L(\E^{1,2}_0)$ is the space of equivalence classes of marked $\E^{1,2}_0$-manifolds 
 of \textbf{linear} holonomy homeomorphic to $\Sigma\times \R$ which are globally hyperbolic,
 Cauchy-maximal, and such that the marking sends $S\times \R$ on $\sing_0$.
\end{defi}

\begin{defi}[Marked $\E^{1,2}_0$-Moduli space] 
Let $\Sigma$ be a compact surface of genus $g$, $S$ be a finite subset of cardinal $s>0$.
 The Moduli space $\M_{g,s}(\E^{1,2}_0)$ is the space of equivalence classes of $\E^{1,2}_0$-manifolds 
  homeomorphic to $\Sigma\times \R$ which are globally hyperbolic,
 Cauchy-maximal, and such that the marking sends $S\times \R$ on $\sing_0$.
\end{defi}

\subsection{First correspondances between Moduli spaces}
 All the Moduli correspondances given in the paper are equivariant under the action of the mapping class group. 
 This is straightforward considering the constructions are explicit, proofs are thus omitted.
 We use brackets, $[M]$, to designate the equivalence class of a manifold $M$. Theorem \ref{barbot} gives a simple correspondance between $\teich_{g,s}$ and $\M_{g,s}^L (\E^{1,2})$ which can be 
given simply. Let $\Sigma$ a surface of genus $g$ and $S$ a finite subset of $\Sigma$ of cardinal $s>0$ with $2g-2+s>0$.
Let $\Gamma$ be the fundamental group of $\Sigma^*$.

\begin{defi}[Suspension] Define the suspension from $\H^2$ to $\E^{1,2}$ 
$$\fonction{\susp_{\H^2}}{\teich_{g,s}}{\M_{g,s}^L (\E^{1,2})}{\left[\Sigma_1^*,\d s^2\right] }{\left[\R_+^*\times \Sigma_1^*,~-\d T^2+T^2 \d s^2 \right]}.$$
 This map comes with a natural embedding of the $\H^2$-surface into its image $\E^{1,2}$-manifold, namely the surface $T=1$.
\end{defi}
The developpement of a point $M$ of $\M_{g,s}^L (\E^{1,2})$ can be chosen such that its image is the future cone $I^+(O)$.
There is a natural time function on $I^+(O)$,  $T:(t,x,y) \mapsto -t^2+x^2+y^2$.

\begin{prop}The map $\susp_{\H^2}$ is bijective of inverse : 

 $$\fonction{\susp_{\H^2}^{-1}}{\M_{g,s}^L (\E^{1,2})}{\teich_{g,s}}{[M]}{\left[(T\circ \D_M)^{-1}(1) \right]} $$
 with $T(t,x,y)=-t^2+x^2+y^2$ and $\D_M$ a developpement of $M$ onto $I^+(O)$. 
\end{prop}
\begin{proof}
 Let $[M]$ be a point of $\M_{g,s}^L (\E^{1,2})$, since the holonomy of $M$ is linear, the developpement of $M$ is the future $I^+(p)$ for some $p\in \E^{1,2}$ which can be chosen 
 to be $O$ the origin of $\E^{1,2}$. The fundamental group $\pi_1(M)$ then acts on $I^+(O)$ via a lattice representation $\rho$
 in $\SO_0(1,2)$. The surface $T=1$ is an isometric embedding of $\H^2$ into $I^+(O)$ thus $ (T\circ \D_M)^{-1}(1)$ is a riemannian 
 surface isometric to $\H^2/\rho$ which is  a point of $\teich_{g,s}$. The suspension of this surface is constructed 
 by identifying its universal cover, namely $\H^2$, with the surface $T=1$ in $I^+(O)$, this shows that the suspension 
 of $\H^2/\rho$ is $I^+(O)/\rho\simeq M$. 
\end{proof}

\begin{rem} For all $[\Sigma_1^*,\d s^2] \in \teich_{g,s}$, the fundamental group of the suspension $\susp_{\H^2}([\Sigma_1^*,\d s^2])$  is canonically 
identified with the one of $\Sigma_1^*$, 
the suspension preserves the marking and 
the holonomy of $\susp_{\H^2}([\Sigma_1^* ,\d s^2])$ is the same as the holonomy of $[\Sigma_1,\d s^2]$.
\end{rem}
We thus obtain the following diagram :
$$\xymatrix{
\teich_{g,s}\ar@<1ex>[rr]^{\susp_{\H^2}}&& \ar@<1ex>[ll]^{\susp^{-1}_{\H^2}} \M_{g,s}^L (\E^{1,2})} $$

\begin{prop} Identifying $T\teich_{g,s}$ and the classes of marked admissible representations $\rho:\Gamma \rightarrow  \isom(\E^{1,2})$,
the holonomy defines an injective map 
$$\Hol : \M_{g,s} (\E^{1,2}) \rightarrow T\teich_{g,s}$$
\end{prop}
\begin{proof}
 Let $[M]\in \M_{g,s}(\E^{1,2})$ and let $\rho$ be its holonomy. $\rho_L$ is admissible by definition of
 $\M_{g,s}(\E^{1,2})$ and by Corollary \ref{cor:admissible} so is $\rho$. 
 The injectivity follows from the construction  
 of the absolutely maximal extension of a globally hyperbolic Cauchy-complete $\E^{1,2}$-manifold which only depends 
 on the holonomy.
\end{proof}

Remains the question of the surjectivity of $\Hol$. From Proposition \ref{prop:abs_ext_2}, it suffices to construct a globally hyperbolic
Cauchy-complete space-time of given admissible marked holonomy. This is the object of section \ref{sec:admissible_space_time}.

In \cite{brunswic_btz_ext}, the maximal BTZ-extension of a $\E^{1,2}$-manifold has been introduced as well as the regular part of 
a $\E^{1,2}_0$-manifold.   Denote these two constructions by $\mathrm{BTZ-ext}$ and $\mathrm{Reg}$. 
    \begin{theo_ext}[Theorem 2 in \cite{brunswic_btz_ext}] \label{theo:BTZ_ext}Let  $A\subset \R_+$, let $M$ be a globally hyperbolic $\E^{1,2}_A$-manifold. 
      
      There exists a maximal BTZ-extension $\overline M$ of $M$. Furthermore it is unique up to isometry.
      \end{theo_ext}

   \begin{theo_ext}[Theorem 3 in \cite{brunswic_btz_ext}] \label{theo:BTZ_Cauchy-completeness}Let $M$ be a globally hyperbolic 
  $\E^{1,2}_{0}$-manifold, the following are equivalent.
  \begin{enumerate}[(i)]
   \item $Reg(M)$ is Cauchy-complete and Cauchy-maximal.
   \item $\mathrm{BTZ-ext}(M)$ is Cauchy-complete and Cauchy-maximal.
  \end{enumerate}
  \end{theo_ext}
 Consider a point $[M]$ of the $\E^{1,2}$-moduli space, its BTZ-extension 
is Cauchy-complete and Cauchy-maximal. In Section \ref{sec:const_BTZ_ext}, we show that $\mathrm{BTZ-ext}(M)$ is Cauchy-compact 
thus a point 
of the $\E^{1,2}_0$-moduli space. Consider now a point $[N]$ of the $\E^{1,2}_0$-moduli space,
its regular part is Cauchy-maximal and Cauchy-maximal. We will show in Section \ref{sec:admissible_hol}
that $\reg(N)$ has admissible holonomy and is thus a point  of  the $\E^{1,2}$-moduli space. 

$$\xymatrix{
 \M_{g,s} (\E^{1,2}) \ar@<1ex>[rr]^{\BTZext}&&\ar@<1ex>[ll]^{\reg} \M_{g,s}(\E^{1,2}_0)
} $$
The constructions $\mathrm{BTZ-ext}$ and $\reg$ being inverse to each other, they will automatically 
define bijections between moduli spaces. In order to complete the picture, we will construct a map 
$\d \susp_{\H^2} : T\teich_{g,s} \rightarrow \M_{g,s}^{\E^{1,2}}$ in Section \ref{sec:admissible_space_time}
Finally, the following diagramq sum-up the situation:
 
 $$\xymatrix{ \teich_{g,s} \ar@<1ex>[rr]^{\susp_{\H^2}}&&\ar@<1ex>[ll]^{\susp_{\H^2}^{-1}}\M_{g,s}^L(\E^{1,2}_0) \ar@<1ex>[rr]^{\mathrm{BTZ-ext}}&&\ar@<1ex>[ll]^{\reg} \M_{g,s}^L(\E^{1,2})}$$

 $$\xymatrix{ T \teich_{g,s} \ar@<1ex>[rr]^{\d\susp_{\H^2}}&&\ar@<1ex>[ll]^{\mathrm{Hol}}\M_{g,s}(\E^{1,2})\ar@<1ex>[rr]^{\mathrm{BTZ-ext}}&&\ar@<1ex>[ll]^{\reg} \M_{g,s}(\E^{1,2}_0)}$$

\section{Spacetime constructions and Moduli spaces}\label{sec:const_spacetime}

\subsection{Maximal BTZ-extension of regular domains}\label{sec:const_BTZ_ext}

We give ourselves a regular domain $\Omega$ invariant under the action of some discrete torsionfree subgroup $G\subset \isom(\E^{1,2})$.
The aim of this section is to give a simple description of the maximal BTZ-extension of $\Omega/G$ as the quotient $\widetilde \Omega/G$ of some 
augmented domain $\widetilde \Omega$. We also prove the Cauchy-compactness of the maximal BTZ-extension of an absolutely 
maximal Cauchy-complete globally hyperbolic $\E^{1,2}$-manifold of admissible holonomy. 
The main result of the section 
are Theorem \ref{theo:max_BTZ_ext}, which gives an explicit construction of the maximal BTZ-extension of a Cauchy-complete globally hyperbolic
$\E^{1,2}$-manifold.
Proposition \ref{prop:fix_line} and Corollary \ref{cor:compact} give precision on Theorem \ref{theo:max_BTZ_ext}
in special cases. 
Also, Definitions \ref{defi:BTZ-line} and \ref{defi:augmented_domain} of augmented regular domain as well as Definition 
\ref{defi:abs_max} of absolutely maximal spacetimes are important in most of what follows.

\subsubsection{Example and augmented regular domain}

Recall the model space of BTZ singularities, namely $\E^{1,2}_0$,  is defined by $\R^3$ with the semi-riemannian metric 
$\d s^2=-\d\tau \d r +\d r^2+r^2 \d \theta$ in cylindrical coordinates. The line $\sing_0(\E^{1,2}_0)=\{r=0\}$ is singular 
and its regular part, $\reg(\E^{1,2}_0)=\{r>0\}$, is a $\E^{1,2}$-manifold.
A developping map $\D : \reg(\E^{1,2}_0)\rightarrow \E^{1,2}$ is given in Proposition 15 and its Corollary in  \cite{brunswic_btz_ext}. 
The image of this developping map (hence every) is the chronological future of some 
lightlike line $\Delta$. Notice that $I^+(\Delta)$ is an open half-space delimited by the lightlike plane $\Delta^\perp$. 
The holonomy is the 
group $\langle \phi \rangle$ where $\phi$ point-wise stabilizes $\Delta$ and  $\D$ induces a homeomorphism 
$$\overline \D :  Reg(\E^{1,2}_0) \xrightarrow{~\sim~} I^+(\Delta)/\langle \gamma \rangle$$ 

These remarks lead to a natural way to construct the maximal BTZ-extension of $\reg(\E^{1,2}_0)$, namely $\E^{1,2}_0$, by quotienting 
$J^+(\Delta)=I^+(\Delta)\cup \Delta$ by $\langle \gamma \rangle$.
The isomorphism $\overline \D$
extends continuously to a bijective map $\overline \D : \E^{1,2}_0 \rightarrow J^+(\Delta)/\langle \gamma\rangle$, by defining 
$\overline \D(\tau,0,0)=(\tau,\tau,0)$.
However, \textbf{ if $J^+(\Delta)$ is endowed with the usual topology, this map is not a homeomorphism}.
This can be seen by taking a sequence of points tending toward $\Delta$ following a horizontal circle intersecting $\Delta$.
The $\tau$ coordinate of the preimage of this sequence goes to $-\infty$. A thinner topology is needed on $J^+(\Delta)$ in order to proceed.

\begin{defi}[BTZ-topology] \label{defi:BTZtop} Let $\Delta$ be a lightlike line in $\E^{1,2}$.
Define the BTZ topology on $J^+(\Delta)$ as the topology generated by the usual open subsets of
	  $J^+(\Delta)$ and the subsets of the form 
	      $$I^+(p)\cup ]p,+\infty[$$   where $p\in \Delta$.
\end{defi}
\begin{prop} \label{prop:quotient_chart}Let $\Delta$ be a lightlike line in $\E^{1,2}$ and $\gamma\in \isom(\E^{1,2})$ fixing $\Delta$ point-wise.
  The map $\overline \D:\BTZ \rightarrow J^+(\Delta)/\langle \phi\rangle$ is a homeomorphism. 
\end{prop}
\begin{proof}
As mentionned before, $\overline \D$ is a bijection. The topology of $\E^{1,2}_0$ is generated by the diamonds 
$\Diamond_p^q:=Int(J^+(p)\cap J^-(q))$ for $p,q\in \E^{1,2}_0$ and the topology of $J^+(\Delta)/\langle \phi\rangle$ is generated by the 
quotient topology 
of  $I^+(\Delta)/\langle \phi\rangle$ and open of the form $\U_p:=(I^+(p)/\langle \phi\rangle)\cup ]p,+\infty[$.

A direct computation gives  $\D^{-1}\left(\U_p\right)= Int(J^+(\D^{-1}(p)))$ which is open. Then, from Corollary 16 
in \cite{brunswic_btz_ext}, $\D$ is continuous. Furthermore, for $p\in \Delta$ and $q\in \BTZ$,
$$\D(\Diamond_{p}^q)=\U_p \setminus \bigcap_{x\in  \Diamond_{p}^q}J^+(\D(x))$$
then $\D(\Diamond_p^q)$ is open. Again from Corollary 16 in \cite{brunswic_btz_ext}, we deduce that $\D$ is open. 
\end{proof}

Let $M$ be Cauchy-complete globally hyperbolic spacetime and let $\Gamma:=\pi_1(M)$. A natural construction of the maximal BTZ-extension of $M$
would then be to consider its developpement $\Omega$ in $\E^{1,2}$ and its holonomy 
$\rho:\Gamma \rightarrow \isom(\E^{1,2})$. 
For $\gamma \in \Gamma$ such that $\rho(\gamma)$ is parabolic, 
 if $\fix(\rho(\gamma))\cap \partial \Omega \neq \emptyset$, we add to $\Omega$ a lightlike ray and quotient out by $\rho$.
 \begin{defi}[BTZ-line associated to parabolic isometry]\label{defi:BTZ-line}
 Let $G$ be a discrete torsionfree subgroup of $\isom(\E^{1,2})$ and $\Omega$ a $G$-invariant regular domain.  
 
 Let $\phi\in G$ parabolic, define the associated BTZ-line $\Delta_\phi$ as the interior of  $\fix(\phi)\cap \partial \Omega$ in $\fix(\phi)$.
  
 \end{defi}
\begin{rem} The BTZ-line associated to a parabolic isometry may be empty. 
 
\end{rem}

\begin{defi}[Augmented regular domain]\label{defi:augmented_domain} Let $G$ be a discrete torsionfree subgroup of $\isom(\E^{1,2})$ and $\Omega$ a $G$-invariant regular domain.  Define

$$\widetilde{\sing_0}(\Omega,G) =  \bigcup _{\phi\in  G ~\mathrm{parabolic}}\Delta_\phi$$
 and $$\widetilde\Omega(G) =\Omega\cup \widetilde{\sing_0}(\Omega,G).$$
 $\widetilde \Omega(G)$ is the augmented regular domain associated to $(\Omega,G)$, we endow it with the BTZ topology. 
\end{defi}
When there is no ambiguity, we shall simply write $\widetilde{\sing_0}$ and $\widetilde \Omega$.
Our aim is now to prove the following theorem.
 \begin{theo} \label{theo:max_BTZ_ext}Let $G$ be a discrete torsionfree subgroup of $\isom(\E^{1,2})$
 and let $\Omega$ be a $G$-invariant regular domain .

 Then $\widetilde\Omega(G)/G$ is endowed with a $\E^{1,2}_{0}$-structure extending the $\E^{1,2}$-structure of $\Omega/G$ and is isomorphic $\mathrm{BTZ-ext}(\Omega/G)$.
\end{theo}
\begin{proof}[Proof for type $(i-iii)$ spacetimes]

 If the group $G$ only consists of spacelike translations. On the one hand $\widetilde \Omega =\Omega$. On the 
 other hand, since the holonomy of a neighborhood of a BTZ-line is parabolic, $\Omega/G$ is BTZ-maximal. The results follows.
 Case $(i)$ is then proved as well as cases $(ii)$ and $(iii)$ with group generated by a spacelike translation.
 
 Assume case $(iii)$ with group generated by a linear parabolic isometry $\phi$. In this case, 
 $\Omega$ is $I^+(\Pi$ for some lightlike plane $\Pi$ and $\Omega/G$ is a annulus 
 $\{(\tau,r,\theta) ~|~ r \in ]R_0,+\infty[\}\subset \E^{1,2}_0$. If $r>0$, then $\Omega/G$ is BTZ-maximal and no
 lightlike line in $\Pi$ is pointwise fixed by $\phi$. Then $\widetilde \Omega=\Omega$ and the results follows. 
 If $r=0$, then there exists a lightlike line $\Delta$ in $\Pi$ point wise fixed by $G$, then $\Omega=I^+(\Delta)$ and 
 $\widetilde \Omega = J^+(\Delta)$. We have $\Omega/G = \reg(\E^{1,2}_0)$ and the results follows from Proposition \ref{prop:quotient_chart}.
 
 Case $(ii)$ with parabolic generator is treated the same way.
 \end{proof}

\subsubsection{Time functions on augmented regular domains}

      We give ourselves a discrete torsionfree isometry subgroup $G\subset \isom(\E^{1,2})$ 
      and a $G$-invariant regular domain $\Omega$. 
      Write $\widetilde \Omega = \widetilde \Omega(G)$, $\widetilde \sing_0 = \widetilde \sing_0(\Omega,G)$ and
      $\sing_0:= \widetilde \sing / G$. Let $M:=\Omega/G$ and $\overline M := \widetilde \Omega/G$.
     \textbf{In this section, we assume $M$ to be a  type $(iv)$ spacetime.}
      
      Let $T:\Omega \rightarrow \R_+^*$ the lift of the Cosmological time function of $M$
      defined in \cite{andersson_time_function}. 
      As mentioned in \cite{MR2499272} Theorem 1.4.1, $T$ is a $\C^1$ Cauchy time function.
      
      \begin{lem} The cosmological $T$ extends continuously to $\widetilde \Omega$ to the map :  
      
      $$\fonction{\widetilde T}{\widetilde \Omega}{\R_+}{p}{\left \{\begin{array}{ll}
                                                                   T(p)& \mathrm{if~}p\in\Omega\\
                                                                   0& \mathrm{if~} p\in \widetilde \sing_0
                                                                  \end{array}
 \right.} $$
     
     \end{lem}
      \begin{proof} 
      Let $p\in \widetilde \sing_0$ and let $\Pi_1$ the lightlike support plane of $\Omega$ at $p$. 
      Since $\Omega/G$ is type $(iv)$, $\Omega$ is in the future of some spacelike plane $\Pi_2$.
       As described in \cite{MR2499272}, for $q\in \Omega$, $T(q)$ is length of the longest past timelike geodesic 
       from $q$. As $q$ goes to $p$, the past timelike geodesics from $q$ in the domain 
       $J^+(\Pi_1)\cap J^+(\Pi_2)\cap J^-(q)$ tends to a segment of lightlike geodesic. Then the length of the longest past timelike 
       geodesic goes to zero and $\lim_{q\rightarrow p} T(q)=0$.
      \end{proof}

      The problem is that $\widetilde T$ is not a time function. Indeed, $\widetilde T$ is non-decreasing for the causal order on $\widetilde \Omega$
      but not increasing. For any $G$-invariant measure $\alpha$ supported on $\widetilde{\sing_0}$ and any positive real $a$, we define a  function
      $$\fonction{T_{\alpha,a}}{ \widetilde{\Omega}}{\R_+^*\cup\{+\infty\}}{p}{\alpha (J^-(p))+a \widetilde T(p)} $$
    Since $\alpha$ is a $G$-invariant measure,  $T_{\alpha,a}$ descends to a function on $M$. Furthermore, if $\alpha$
    is absolutely continuous with respect to Lebesgue measure on $\widetilde \sing_0$ and $T_{\alpha,a}<+\infty$ then 
    $T_{\alpha,a}$ is increasing thus a time function. If moreover, $\alpha(\Delta)=+\infty$ for every BTZ-line $\Delta$
    then $T_{\alpha,a}$ will be a Cauchy-time function.
      One could choose the Lebesgue measure on some BTZ-line  $\Delta$  and then sums all the 
      translation by $G/\stab(\Delta)$, however this gives an infinite $T_{\alpha,a}$ in general. 
      We thus consider such measures with cut-off below a certain point and the sum their translations.
      Lemmas \ref{cor:discrete} and \ref{lem:stab} will ensure the sum induces a well defined, finite $T_{\alpha,a}$.

\begin{defi}
      Let $\Delta$ be a line of $\E^{1,2}$. Define $G_\Delta$ the set of elements of $G$ stabilising $\Delta$
      point-wise.
      $$ G_\Delta := \bigcap_{p\in \Delta} \stab_G(p)$$
\end{defi}

\begin{lem}\label{lem:stab} Let $\Delta$ be a lightlike line such that $G_\Delta \neq \{1\}$.
      Then $G_\Delta$ only contains parabolic isometries and an element of $G$ stabilizing $\Delta$ set-wise is in $G_\Delta$.
\end{lem}
\begin{proof}
      
      To begin with, since $\Omega/G$ is type $(iv)$, then, 
      by Theorem \ref{theo:barbot_mess}, $G$ and $L(G)$ are discrete and $L_{|G}$ is injective. 
      An isometry which stabilises $\Delta$ point wise is conjugated to a linear isometry fixing point-wise a lightlike line. 
      It is thus conjugated to a linear parabolic isometry.
      Then, $G_{\Delta}$ is a non trivial discrete subgroup of $\bigcap_{p\in\Delta}\stab_{\mathrm{Isom}(\E^{1,2})}(p)\simeq \R$,
      then  $G_\Delta$ is monogene. Let $\phi_\Delta$ be a generator of  $G_\Delta$.

      Let $\phi\in G\setminus G_\Delta$  be such that $\phi \cdot \Delta =\Delta$. Let $P\leq Q$ in $\Delta$ be such that 
      $\phi P = Q$.
      Up to conjugating by a translation, we can assume $P=O$ the origin of $\E^{1,2}$ so that $\phi_L P=P$, 
      $\tau_\phi=\overrightarrow {PQ}$ and $\phi_\Delta$ is linear.
      Since $\phi \Delta =\Delta$, $\phi_L \overrightarrow \Delta = \overrightarrow \Delta$ and thus $\overrightarrow \Delta$
      is a lightlike eigen direction of $\phi_L$. Then $\phi$ is parabolic or hyperbolic.
      The group generated by $\phi_L$ and $\phi_\Delta$ is a discrete subgroup of $\SO_0(1,2)$ fixing a point of the boundary 
      of $\H^2$, it is thus monogeneous and let $\psi$ be a generator. 
      There exists $p,q\in \Z$ such that $\psi^p= \phi_\Delta$ and $\psi^q=\phi_L$. Since $\phi_\Delta$ is parabolic,
      so is $\psi$ and thus so is $\phi_L$. We have $\phi_\Delta^{q}\phi^{p}=\tau_\phi$ then $\tau_\phi \in G$.
      Since $L\circ \rho$ is faithful and $L(\tau_\phi)=0$ we have $\tau_\phi=0$, then $\phi=\phi_L$ and thus $\phi \in G_\Delta$.

      \end{proof}
\begin{cor} \label{cor:stab} Let $\Delta$ be a BTZ-line of $\widetilde \Omega$ and let $\psi\in G$. 
If there exists $p\in \Delta$ such that $\psi p \in \Delta$ and $\psi p \neq p$ then $\psi \in G_\Delta$.
\end{cor}
\begin{proof} 
 Let $q\in \Delta$, $\psi q = \psi(q-p+p)=\psi_L(q-p)+\psi p$. Since $\psi p$ and $\psi q$  are in $\sing_0$,  
 either $\psi p - \psi q$ is spacelike or $p,q$ belongs to the same BTZ-line. 
 The former is not possiblie since $q-p$ is lightlike, therefore $\psi q \in \Delta$. Then, $\psi$ stabilises $\Delta$ set-wise 
 and from Lemma \ref{lem:stab}, $\psi \in G_\Delta$.
\end{proof}

\begin{lem}\label{lem:discrete_effective}Let $\Delta$ be a BTZ-line of $\widetilde \Omega$.  For all $p \in \Delta$, 
there exists $\lambda>0$ such that 
$$\forall q\in \widetilde \Omega,\quad \# \left(Gp \cap J^-(q) \right) \leq (1+\lambda \widetilde T(q))^2 $$
\end{lem}
\begin{proof}
Let $q\in \widetilde \Omega$, if $q\in \widetilde \sing_0$ by Lemma \ref{lem:stab}, $\#\left(Gp \cap J^-(q)\right) \leq 1$.

Let $p\in \Delta$, let $p_*=\inf(\Delta)$ and let $u=p-p_*$. The vector $u$ is future lightlike vector and 
$\Delta = p_* + \R_+^* u$.  
For $v$ lightlike, define $h_{v}=J^+(v)\cap \H^2$. 
The set $\{h_{tv} : t>0\}$ is exactly the set of horocyles centered at $v$. Since $L(G)$ is discrete, $\H^2/L(G)$ is
a complete $\H^2$-manifold and there exists an embedded horocycle around the cusp associated to $u$. 
Let $\lambda >0$ such that $h_{\lambda u}$ is embedded. Let $\overrightarrow n$ be the vertical future unit timelike vector.
Let $\phi \in G$
\begin{eqnarray} 
 \phi p \in J^-(q) &\Leftrightarrow & \phi p_* + \phi_L u \in J^-(q)\\
 &\Leftrightarrow &  \phi_L u \in J^-(q-\phi p_*) \\
 &\Rightarrow&\left| \langle \phi_L u | \overrightarrow n\rangle\right| \leq \left| \langle q-\phi p_* | \overrightarrow n\rangle   \right|\\
 &\Rightarrow& \left| \langle\phi_L u|\overrightarrow n  \rangle \right| \leq T(q)
\end{eqnarray}

On the one hand, for $v$ lightlike,
the stereographic projection of $\H^2$ onto the Poincaré disc on $\overrightarrow n^\perp$  sends an horocycle 
$J^+(v)\cap \H^2$ to a Euclidean circle of radius $(1+|\langle v |\overrightarrow n \rangle|)^{-1}$.
 On the other hand, the horocycles $h_{\lambda \phi_L u}$ are disjoint for $\phi_L \in L(G)/L(G_\Delta)$. 
If $\phi p \in J^-(q)$, then the radius of $h_{\lambda \phi_L u}$ is greater than 
$(1+\lambda T(q))^{-1}$. Since the total area of the disjoint horoball is less than $\pi$, 
there exists at most $(1+\lambda T(q))^2$ such $\phi_L \in L(G)/L(G_\Delta)$. 
Since $\Omega/G$ is  type $(iv)$, then $L_{|G}$ is injective  and the result follows.

\end{proof}

\begin{cor}\label{cor:discrete} Let $\phi \in G$ parabolic. 
Then for all $p \in \Delta_\phi,~ G p$ is discrete.
\end{cor}

        \begin{prop} There exists a measure  $\alpha$  $\widetilde{\sing_0}$ such that for all $a\in \R_+^*$,  
        \label{prop:time_function}
        \begin{itemize}
         \item $T_{\alpha,a}$ is $\C^1$ on $\Omega$ and $\C^0$ on $\widetilde \Omega$;
         \item $T_{\alpha,a}$ is $G$-invariant.
         \item $T_{\alpha,a}$ is a Cauchy time function on $\widetilde \Omega$
         
         \end{itemize}
    \end{prop}
 \begin{proof}
 Choose a set of representative $(\Delta_i)_{i\in I}$ of $\widetilde \sing_0$. The set $I$ is countable we can thus assume $I\subset \N$ and for each $i\in I$,
 choose a decreasing sequence $(p_n^{(i)})_{n\in \N} \in \Delta_i^\N$ such that
 $\lim_{n\rightarrow +\infty}p_n^{(i)} = \min( \Delta_i)$. 
 Let $N(i,n)$ be the number of triplet  $(j,k,\psi)$ with $j \leq i$ and $k\leq n$   and $\psi \in G/G_{\Delta_j}$ such
 that $\psi p_k^{j} \in $
 We can choose a family $(\varphi_n^{(i)})_{n\in \N,i\in I}$ such that for all $n\in \N$ and $i\in I$, 
 
 \begin{enumerate}[$(i)$]
  \item $\varphi_n^{(i)}$ is in $\C^1(\Delta_i, \R_+)$, 
  \item $\|\varphi_n^{(i)}\|_{\C^1} \leq 1$
  \item $\lim_{x\rightarrow +\infty} \varphi_n^{(i)}(x)=1$
  \item $\forall x\in \Delta_i, \varphi_n^{(i)}(x)=0 \Leftrightarrow x\leq p_n^{(i)}$
 \end{enumerate}
 
 Choose a geodesic parametrisation of $\Delta_i$ for each $i\in I$ and let $\lambda_i$ the image of the Lebesgue measure 
 on $\R_+^*$ by this parametrisation. From Lemma \ref{lem:discrete_effective}, 
 for $n\in \N$ and $i\in I$, let $\mu_n^{(i)}\geq 1$ be such that 
 $$\forall q\in \Omega, \quad\# \left\{\phi \in G/G_\Delta ~|~ \phi p_n^{(i)} \in J^-(q) \right\} \leq (1+\mu_n^{(i)} T(q))^2 $$
 Let $$\alpha =\sum_{i\in I }\sum_{\psi \in G/G_{\Delta_i}}\sum_{n\in \N} \omega_n^{(i)} \psi\#(\varphi_n^{(i)} \lambda_i) $$
 where  $\omega_n^{(i)} = \frac{2^{-i-n}}{\lambda_i\left(J^-\left(p_0^{(i)}\right) \right) \mu_n^{(i)}}$.

 Define for $i\in I$ and $n\in \N$, $$\fonction{\alpha_n^{(i)}}{\widetilde \Omega}{\R_+}{p}{\sum_{\psi \in G/G_\Delta}} \varphi_n^{(i)} \lambda_i\left(J^-(\psi p)\right)$$
 The sum is locally finite thus $\alpha_n^{(i)}$ is $\C^1$ and finite. Furthermore, 
 for all $q\in \widetilde \Omega$ : 
 $$\left\|\alpha_n^{(i)} \right\|_{\C^1(J^-(q))}\leq \lambda_i\left(J^-\left(p_0^{(i)}\right) \right) \left(1+\mu_n^{(i)}T(q)\right)^2.$$
 Then, for all $q\in \widetilde\Omega$ :
$$
\sum_{i\in I} \sum_{n\in \N} \left\|\omega_n^{(i)}\alpha_n^{(i)}\right\|_{\C^1(J^-(q))}\leq \sum_{i\in I} \sum_{n\in \N}  2^{-i-n} (1+T(q))^2=4 (1+T(q))^2.
 $$
 Thus, the sum $\sum_{i\in I} \sum_{n\in \N} \omega_n^{(i)}\alpha_n^{(i)}$ is normally convergent on compact subset of $\widetilde \Omega$
for the $\C^1$ norm and is thus $\C^1$.

It remains to prove $T_{\alpha,a}$ is Cauchy, i.e that $T_{\alpha,a}$ is surjective and increasing on inextendible causal curves. 
Let $c:\R\rightarrow \widetilde \Omega$ be an inextendible causal curve, define $\Delta=c\cap \widetilde \sing_0$ and $c^0 = c\cap \Omega$.
The two pieces $\Delta$ and $c^0$ are connected and $\Delta$ is in the past of $c^0$. The function 
$T$ is increasing on $c^0$ then so is $T_{\alpha,a}$. Since $\alpha$ is absolutely continuous with respect to Lebesgue measure 
on $\widetilde\sing_0$ then $ T_{\alpha,a}$ is increasing on $\Delta$.
When $t\rightarrow -\infty$, $\widetilde T(c(t))$ and $\alpha(J^-(c(t)))$ go to $O$ thus $T_{\alpha,a}(c(t))$ goes to 0.
If $c^0=\emptyset$, then $\bigcup_{t>0} J^-(c(t))$ is a connected component component of $\widetilde \sing_0$ and by condition 
$(iii)$, $\alpha(J^-(c't))$ goes to $+\infty$. If $c^0\neq \emptyset$,  then for $c^0$ is a non-empty inextendible future causal curve
of $\Omega$ and thus $\lim_{t\rightarrow +\infty} T(c(t)) = +\infty$. 
In any case, $\lim_{t\rightarrow +\infty } T_{\alpha,a}(c(t))=+\infty $.
Finally, $T_{\alpha,a}$ is a Cauchy time function on $\widetilde \Omega$.
\end{proof}

\subsubsection{Construction of the maximal BTZ-extension  of a regular domain}

     We give ourselves a discrete torsionfree isometry subgroup $G\subset \isom(\E^{1,2})$ and
     a $G$-invariant regular domain $\Omega$. We assume $\Omega/G$ is of type $(iv)$. 
     Write (see Definition \ref{defi:augmented_domain}) $\widetilde \Omega = \widetilde \Omega(G)$, 
     $\widetilde \sing_0 = \widetilde \sing_0(\Omega,G)$.
     Let $M:=\Omega/G$, $\MM := \widetilde \Omega/G$ with the quotient topology and let $\pi : \widetilde \Omega\rightarrow \MM$ be the natural projection.

\begin{prop}\label{prop:Hausdorff}
  $\MM$ is Hausdorff.
\end{prop}
\begin{proof}
	Let $\alpha$ be a measure given by Proposition \ref{prop:time_function} and let $a\in \R_+^*$. 
	Let $p,q\in \widetilde \Omega$ be such that $\pi(p)\neq  \pi(q)$. 
	\begin{itemize}
	  \item If $p$ and $q$ are in $\Omega$, then $\pi(p)$ and $\pi(q)$ are in  $M$ and since $M$ is Hausdorff,
	  then $\pi(p)$ and $\pi(q)$ are separated in $M$. And thus in $\pi(p)$ and $\pi(q)$ are separated in $\MM$. 
	  \item If $p$ and $q$ are in $\widetilde \sing_0$, then either $T_{\alpha,a}$ separates $p$ and $q$ or $p,q$ are not 
	  on the same BTZ-line. In the latter case we can multiply $\alpha$ by a $G$-invarariant function equal to $1/2$ on the orbit of 
	  the BTZ-line of $p$ and $1$ on the other BTZ-lines to obtain $T_{\alpha',a}(p)\neq T_{\alpha',a}(q)$. 
	  Since $T_{\alpha,a}$ is $G$-invariant and continuous, $\pi(p)$ and $\pi(q)$ are separated.
	  \item If $p\in \Omega$ and $q\in \widetilde \sing$, then one can change the $a$ parameter to obtain $T_{\alpha,a}(p)\neq T_{\alpha,a}(q)$. Then  we can use the same argument as before.
	\end{itemize}   
\end{proof}

\begin{prop}\label{prop:structure}
  $\MM$ has  a $\E^{1,2}_{0}$-structure such that $Reg(\MM)=M$. Furthermore, this structure is globally hyperbolic.
\end{prop}
\begin{proof}
	The $\E^{1,2}$-atlas of $M$ gives an atlas of $\MM \setminus (\widetilde\sing_0)/G$, 
	it suffices to construct charts around points 	of $\widetilde \sing_0 /G$. 
	Beware that hereafter, the topology on 
	$\widetilde \Omega$ is the BTZ-topology given in Definition \ref{defi:BTZtop}.	
	Let $\Delta$ be a connected component of $ \widetilde \sing_0$ (i.e. a BTZ-line of $\widetilde \Omega$) and $p\in \Delta$ be some point on it. 
	Let $G^*:=G\setminus G_\Delta$.
	We can choose the map $\overline \D$ of Proposition \ref{prop:quotient_chart} such that $\Delta\subset \overline \D(\sing_0(\E^{1,2}_0))$.
	Denote $q=\overline \D^{-1}(p)$.
	We now  construct a $G_\Delta$-invariant  neighborhood $\W$ 
	of $p$, open for the BTZ-topology and disjoint from its image by $G^*$. 
	Then $\overline \D$ will induces a homeomorphism between  some open neighborhood 
	of $q=\overline \D^{-1}(p)$ and $\W/G_\Delta$.
	
	      \begin{figure}\label{fig:dvp_neigh}
    
           \includegraphics[width=4.9cm]{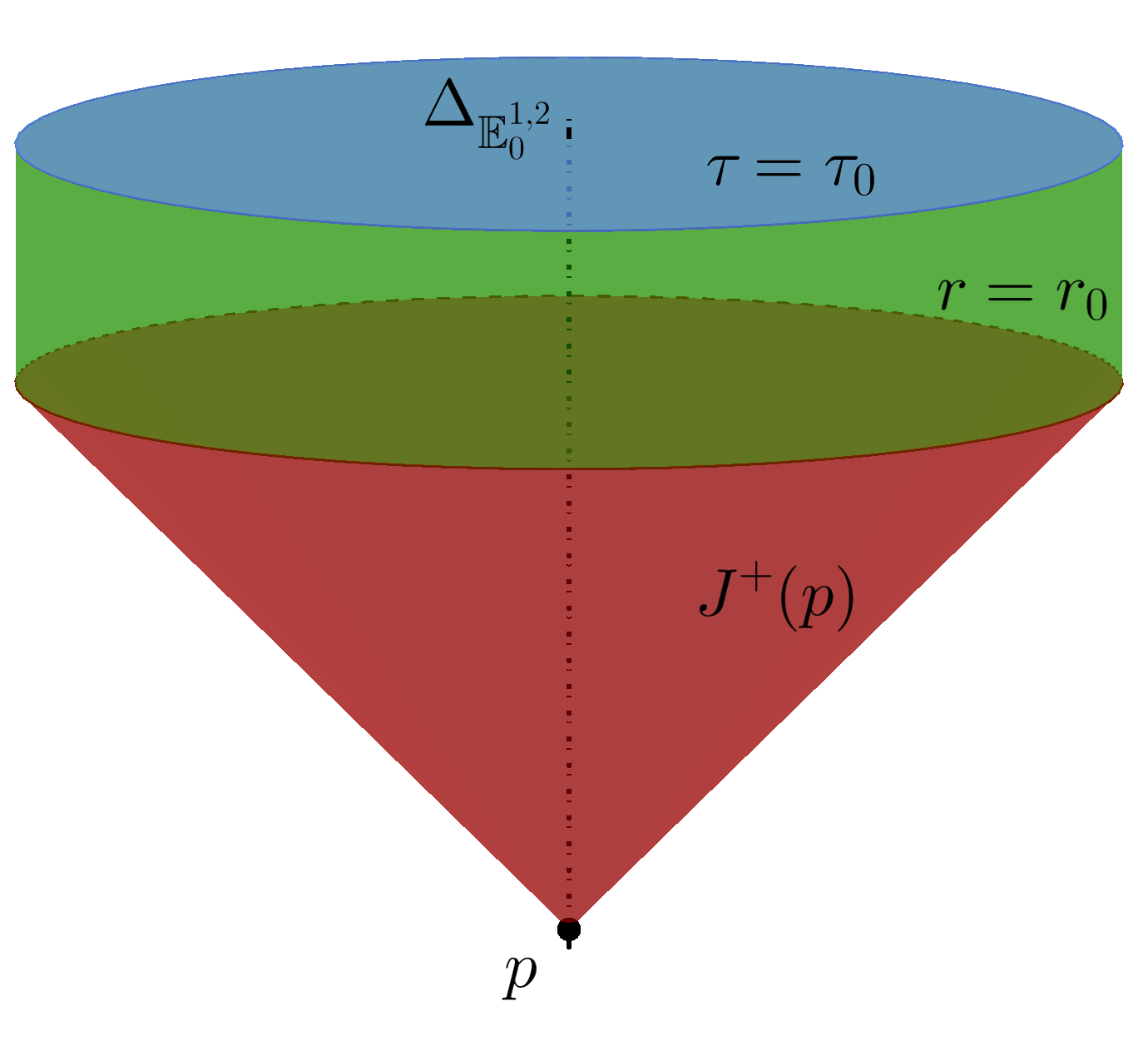}
        \includegraphics[width=6.9cm]{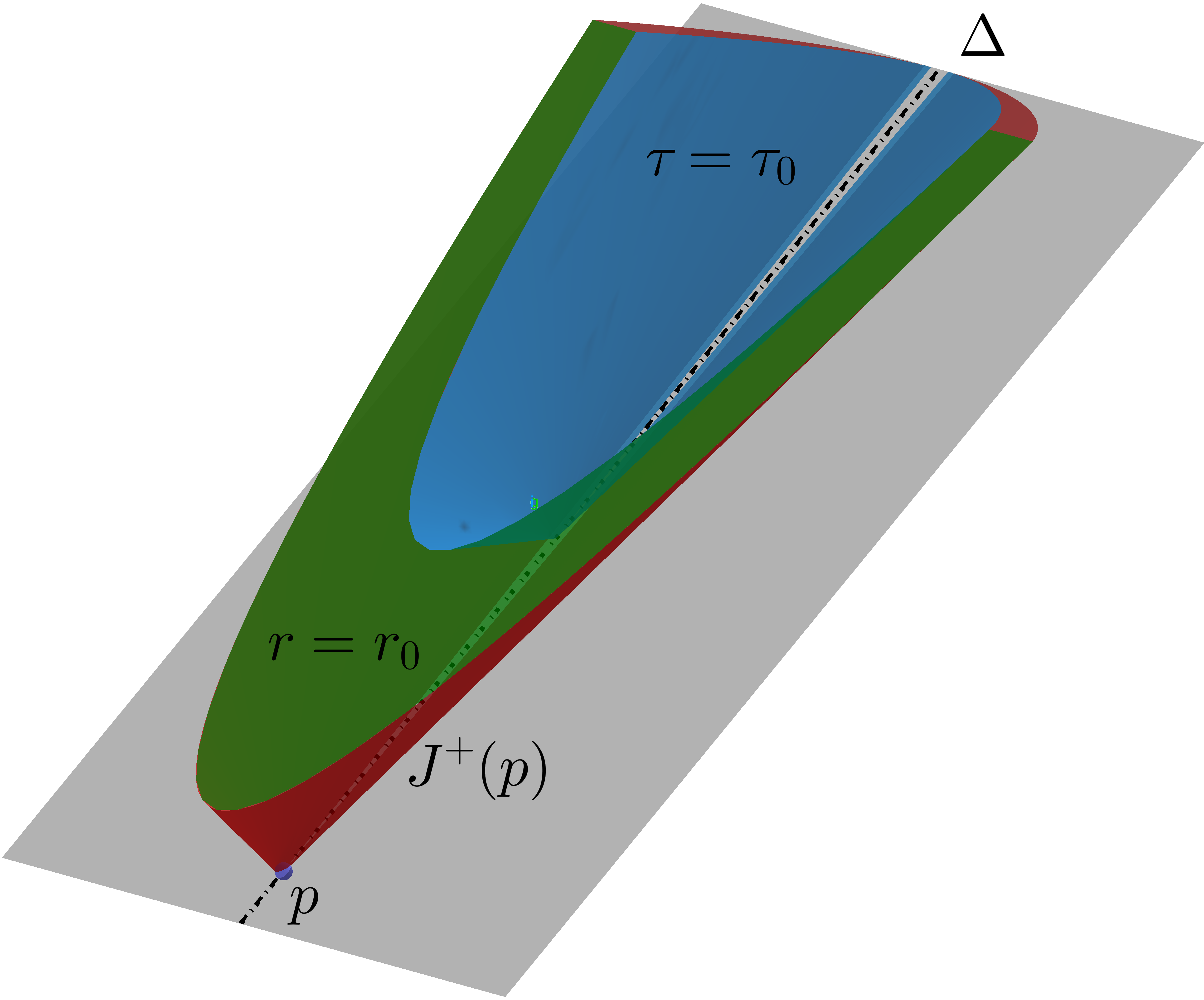}
        \begin{caption}{Tubular neighborhood of a BTZ point and its development} On the left, a tubular subset of $\E^{1,2}_0$. On the right its development into $\E^{1,2}$. Colors are associated 
        to remarkable sub-surfaces and their developments.\label{fig:neighborhood_BTZ}
        \end{caption}
    \end{figure}
 
	 Let $q_* \in J^-(q)$ and denote $p_*=\overline D(q)$ so that $p_*\in J^-(p)$ and let
	 $$\V=Int\left(J^+(q_*)\cap\{r<R, \tau<\tau^* \}\right) \quad \mathrm{and}\quad \U = \pi^{-1}\left(\overline\D(\V)\right)$$
	  where $(\tau,r,\theta)$ are the  cylindrical coordinates of $\E^{1,2}_0$ and  $R>0$ and  $\tau^*>\frac{1}{2}R$. 
	 See Figure \ref{fig:dvp_neigh} below which depicts $\V$ and $\U$.
	Now, choose $\U_0\subset \U$ some relatively compact open domain such that $G_\Delta \cdot \U_0=\U$, then: 
	
	\begin{eqnarray}
	 \U \setminus G^*\U&=&\U \setminus \left(G^*\U\cap \U\right) \\
	 &=& (G_\Delta \U_0) \setminus  G_\Delta\left((G^*  \U) \cap \U_0\right)\\
	 &=& G_\Delta  \left(\U_0\setminus (G^* \U)\cap \U_0\right)
	\end{eqnarray}
	Since $\U\subset J^+(p_*)$, for $\phi\in G$,
	if $(\phi\U)\cap \U_0\neq\emptyset $ then $\phi p_* \in J^-(\U_0)\subset J^-(\overline{\U_0})$.
	From Corollary \ref{cor:discrete} and since $J^-(\overline \U_0)$ is compact,
	the set  $G p_* \cap J^-(\overline \U_0)$ is finite, say $\{\phi_0p_*,\phi_1 p_*,\cdots,\phi_np_*\}$ with $\phi_0=1$.
	Therefore : 
\begin{eqnarray}
	 \U \setminus G^*\U&=&G_\Delta \cdot \left(\U_0\setminus (\U_0\cap G^* \U)\right)\\
	 &=& G_\Delta\cdot\left(\U_0\setminus \bigcup_{i=1}^{n} \phi_i \U\right)\\
	 &\supset& G_\Delta\cdot\left(\U_0 \setminus \bigcup_{i=1}^n J^+(\phi_i p_*)\right)
	\end{eqnarray}
	Notice that from Lemma \ref{lem:stab}, for all $\phi\in G$ such that $\phi p_* \in J^-(p_*)$, 
	we have $\phi \in G_\Delta$. 
	Thus the only $\phi_i$ 	such that $\phi_i p_*\in J^-(\overline \U_0)$ is $\phi_0=1$.
	Let $\W:=G_\Delta\left(\U_0 \setminus \bigcup_{i=1}^n J^+(\phi_i p_*)\right)$, 
	then $\W$  is an open subset of $\U$, disjoint from $G^* \W$  and containing $p$. 

	A function $T_{\alpha,a}$ from Proposition \ref{prop:time_function} is a $G$-invariant Cauchy-function of $\widetilde\Omega$.
	It thus induces a Cauchy function of $\MM$ which proves that $\MM$ is globally hyperbolic.
\end{proof}

 \setcounter{theo}{0}
 \begin{theo} Let $G$ be a discrete torsionfree subgroup of $\isom(\E^{1,2})$
 and let $\Omega$ be a $G$-invariant regular domain .

 Then $\widetilde\Omega(G)/G$ is endowed with a $\E^{1,2}_{0}$-structure extending the $\E^{1,2}$-structure of $\Omega/G$ and is isomorphic $\mathrm{BTZ-ext}(\Omega/G)$.
\end{theo}

\begin{proof}[Proof for type $(iv)$ spacetimes]
  Propositions \ref{prop:Hausdorff} and \ref{prop:structure} prove that $\Omega/G \rightarrow \widetilde \Omega/ G$ is 
  a BTZ-embedding. It remains to prove that $\widetilde \Omega/G$ is BTZ-maximal.
  Consider the maximal BTZ-embedding $\widetilde \Omega/G \xrightarrow{i} N$, 
  take a point $p\in\sing_0(N)$, a compact diamond neighborhood $\U$ around $p$ included in a chart around $p$ 
  and some loop $c:=\{r=R_0,\tau=\tau_0\}$   around the line $\sing_0(\U)$. Let $\D$ be the developping map of $\reg(N)$.
  The image of the holonomy of $\reg(\U)$ is generated
  by a parabolic isometry $\phi$ fixing a lightlike line $\Delta$ which intersects the boundary of $\D(\widetilde{\reg(\U)})$
  along a segment $[p_*,p^*]$.  Let $\V:=\D(\widetilde{\reg(\U)})$,
  we can assume that $(G\setminus G_\Delta) \V \cap \V = \emptyset$. 
  By Proposition \ref{prop:quotient_chart}, $\D$ induces a homeomorphism 
  $\overline \D : \U \rightarrow (\V\cup[p_*,p^*])/G_\Delta $ and 
  $i_{|\V/G_\Delta}$ is a continuous section of $\overline \D$ on $\V/G_\Delta$. Therefore, by continuity, 
  $\overline \D^{-1}=i_{|(\V\cup [p_*,p^*])/G_\Delta}$ and $p$ is in the image of $i$. 
  Finally, $i$ is surjective.
\end{proof}

\subsubsection{Absolutely maximal singular spacetimes}
We introduce the notion of absolutely maximal $\E^{1,2}_A$-manifolds which will prove relevant in our description of Lorentzian 
Moduli spaces. See Definition 21 in \cite{brunswic_btz_ext} for a definition of $\E^{1,2}_A$-manifolds.

\begin{defi}\label{defi:abs_max} A connected Cauchy-complete globally hyperbolic $\E^{1,2}_A$-manifold $M_1$ is absolutely-maximal if all $\E^{1,2}_A$-embedding 
$M_1\rightarrow M_2$, with $M_2$ Cauchy-complete globally hyperbolic and connected, is onto.
 
\end{defi}
\begin{rem} Beware that the absolute maximality depends strongly on the category of $\E^{1,2}_A$-manifold you consider. 
Indeed, an absolutely maximal $\E^{1,2}$-manifold may have a BTZ-extension and thus not absolutely as  
$\E^{1,2}_0$-manifold. For instance :  $\reg(\E^{1,2}_0)$ is absolutely maximal as $\E^{1,2}$-manifold but embeds into $\E^{1,2}_0$, thus $\reg(\E^{1,2}_0)$ is 
not absolutely maximal as $\E^{1,2}_0$-manifold.
\end{rem}

 We do not claim that there exists a unique  absolutely maximal extension for every $\E^{1,2}_A$-manifold.
 A theorem of existence and unicity is known for Cauchy-complete $\E^{1,2}$-manifold as a consequence of Theorem \ref{theo:barbot_mess}.
 We now use Theorem \ref{theo:max_BTZ_ext} 
 to extend it for Cauchy-complete $\E^{1,2}_0$-manifolds.

\begin{prop}\label{prop:abs_max_fund} Let $M$ be a connected Cauchy-complete globally hyperbolic $\BTZ$-manifold then there exists a 
$\BTZ$-manifold $\overline M$ absolutely maximal in which $M$ embeddeds. Moreover, $\overline M$ is unique up to isomorphism. 
\end{prop}
\begin{proof}
 Let $M_1$ be the maximal Cauchy-extension of $M$, it exists and is unique (see for instance \cite{Particles_1}).
 From Theorem \ref{theo:BTZ_Cauchy-completeness},  $\reg(M_1)$ is a  Cauchy-complete and Cauchy-maximal globally 
 hyperbolic $\E^{1,2}$-manifold. Let $M_2$ be the maximal BTZ-extension of the absolutely maximal extension of $\reg(M_1)$ among $\E^{1,2}$-manifolds.
 Let $\widetilde \Omega_1$ (resp. $\widetilde \Omega_2$ ) be the augmented regular domain associated to 
 $\reg(M_1)$ (resp. $\reg(M_2)$). 
 These augmented regular domain can be chosen such that $\widetilde \Omega_1 \subset \widetilde \Omega_2$. 
 From Theorem \ref{theo:max_BTZ_ext}, this inclusion induces an embedding $i:M_1 \rightarrow M_2$. 
 Let $M_3$ be a connected Cauchy-complete Cauchy-maximal $\E^{1,2}_0$-manifold and an embedding $j:M\rightarrow M_3$.
 We can extend $j$ to an embedding $M_1 \rightarrow M_3$. 
 We thus obtain an embedding $\reg(M_1) \rightarrow \reg(M_3)$ and thus by absolute maximality of $\reg(M_2)$, an embedding 
 $\reg(M_3)\rightarrow \reg(M_2)$.
 Then we obtain an inclusion 
 $\widetilde \Omega_3 \subset \widetilde \Omega_2$ with $\widetilde \Omega_3$ the augmented regular domain associated to $\reg(M_3)$.
 By Theorem \ref{theo:max_BTZ_ext}, we then get an emdedding $M_3 \rightarrow M_2$ such that the following diagram commutes
 $$\xymatrix{ M\ar[d]^j\ar[r]&M_1\ar[d]^i \\ M_3 \ar[r]& M_2 }$$
 
 If we have an embedding $f:M_2 \rightarrow M_3$ then $M$ embedds into $M_3$ and thus we obtain a map $g:M_3\rightarrow  M_2$. 
 The commutative diagram then implies that $g\circ f= Id_{M_2}$ and thus $f$ is surjective. Thus $M_2$ is absolutely maximal.
 If $M_3$ is absolutely maximal then the map  $M_3 \rightarrow M_2$ is surjective thus an isomorphism.   Thus $M_2$ is unique 
 up to isomorphism.

\end{proof}

\begin{lem}\label{lem:vector} Let $M$ be a $\E^{1,2}_A$ manifold, there exists a vector field on $M$ 
such that 
\begin{itemize}
 \item $X$ is $\C^1$ and non singular;
 \item $X$ is future causal;
 \item For $p\in \sing(M)$, $X_p$ is parallel to the direction of the singular line through $p$.
\end{itemize}
\end{lem}
\begin{proof}[Sketch of proof] 
Let $\V=\reg(M)$, $\V$ has a time order relation.
Thus, from Lemma 32 p145 of \cite{oneil}, there exists a timelike hence causal and non-singular $\C^1$ vector field $X_\V$ on $\V$. 
Take a family of disjoint  chart neighborhoods $(\U_i)_{i\in I}$ of the singular lines $(\Delta_i)_{i \in I}$.
On each  $\U_i$, define $X_i$ a constant vector field parallel to $\Delta_i$ then construct a partition of unity $(\varphi_{\U_i}: i\in I, \varphi_\V)$ associated 
to  the open cover $(\U_i : i\in I, \V)$ such that $\varphi_{\U_i}=1$ on a neighborhood of $\Delta_i$. The vector field $\sum_{i\in I} \varphi_{\U_i} X_i+\varphi_{\V} X_\V$ 
statifies the wanted properties.
\end{proof}

\begin{lem}\label{lem:cauchy_max_abs_max} Let $M$ be a globally hyperbolic $\E^{1,2}_A$-manifold. If $M$ is Cauchy-compact, then the following are equivalent \label{prop:c-compact_maximal} 
\begin{enumerate}[(i)]
 \item $M$ is Cauchy-maximal;
 \item $M$ is absolutely maximal.
\end{enumerate}
\end{lem}
\begin{proof} Assume $M$ is Cauchy-compact.
 \begin{itemize}
  \item Assume $M$ is absolutely maximal.  Since a Cauchy-embedding is an embedding of $\E^{1,2}_A$-manifold, in particular every 
  Cauchy-embedding is onto and thus $M$ is Cauchy-maximal.
  \item Assume $M$ is Cauchy-maximal and Consider $M\xrightarrow{i} M'$ a $\E^{1,2}_A$-embedding with $M'$ globally hyperbolic.
  Take $\Sigma_1$ (resp. $\Sigma_2$) a spacelike Cauchy-surface of $M_1$ (resp. $M_2$), such a surface exists from Theorem 1 in 
  \cite{brunswic_btz_ext}. 
  Take $X$ a vector field on $M'$ given
  by Lemma \ref{lem:vector}. For every $p\in \Sigma_1$, the line of flow of $X$ through $p$ intersects $\Sigma_2$ exactly once. 
  The map $f:\Sigma_1\rightarrow \Sigma_2$ defined this way is a local homeomorphism. Since $\Sigma_1$ is compact, then $f$
  is proper, $f$ is thus a covering and thus $f^* : \pi_1(\Sigma_1) \rightarrow \pi_1(\Sigma_2)$ is surjective. Finally,
  $i^*:\pi_1(M_1) \rightarrow \pi_1(M_2)$ is onto. 
  The proof of Lemma 45 p427 of \cite{oneil} applies to the context of $\E^{1,2}_A$-manifolds 
  and from this Lemma we deduce  that $\Sigma_1$ is achronal, hence acausal.
  
  From Lemma 43 p426 of \cite{oneil}, the Cauchy developpement of $\Sigma_1$ is open and, 
  since $\Sigma_1$  is compact, it is also closed. By connectedness of $M_2$, the Cauchy developpement 
  of $\Sigma_1$ is the whole $M_2$ thus $i$ is a Cauchy-embedding. However, $M$ is Cauchy-maximal thus $i$ is surjective.
  
 \end{itemize}

\end{proof}

\begin{rem} Neither the Cauchy-compacity nor the Cauchy-maximality of a spacetime $M$ depend on the category of $\E^{1,2}_A$-manifold 
in which we are considering $M$. Therefore, a Cauchy-compact and Cauchy-maximal $\E^{1,2}_A$-manifold is absolutely maximal whatever
the category in which it is considered.
\end{rem}

\begin{lem} \label{lem:abs_max_compact}
Let $N$ be a Cauchy-complete Cauchy-maximal globally-hyperbolic $\E^{1,2}$-manifold and let $\overline N$ be 
its absolutely maximal extension among $\E^{1,2}_0$-manifold.

If $\overline N$ is Cauchy-compact, then the following are equivalent: 
\begin{enumerate}[(i)]
 \item $N$ is absolutely maximal among $\E^{1,2}$-manifolds;
 \item $\mathrm{BTZ-ext}(N)=\overline N$
 \item $\mathrm{BTZ-ext}(N)$ contains a point of each BTZ-line of $\overline N$.
\end{enumerate}
\end{lem}
\begin{proof}
Let $i: N \rightarrow \overline N$ be an embedding, $i(N)\subset \reg(\overline N)$. 
\begin{itemize}
 \item If $N$ is absolutely maximal, then  $i(N)=\reg(\overline N)$ and $\mathrm{BTZ-ext}(N)\simeq\overline N$. Thus $(i)\Rightarrow (ii)$.
  \item  Assume $\mathrm{BTZ-ext}(N)\simeq\overline N$, we have $\reg(\overline N)=i(N)$.
  Let $j:N\rightarrow N'$ be the absolutely maximal extension of $N$ among 
  $\E^{1,2}$-manifold. From Proposition \ref{prop:abs_max_fund}, the absolutely maximal extension of $N'$ 
  is $\overline N$ and we have an embedding $k:N'\rightarrow \overline N$. 
  Since, $k(N') \subset \reg(\overline N)=i(N)$, then $\phi := i^{-1}\circ k\circ j$  is an automorphism  
  of  $N$.  Finally,  $j=k^{-1}\circ i \circ \phi$ is surjective 
and thud $N=\reg(\overline N)$ is absolutely maximal.
This proves $(i)\Leftarrow (ii)$.

\item  
Assume $\mathrm{BTZ-ext}(N)=\overline N$, then trivially $\mathrm{BTZ-ext}(N)$ contains a point 
of each BTZ-line of $\overline N$. Therefore, $(ii)\Rightarrow (iii)$.

\item Assume $\mathrm{BTZ-ext}(N)$ contains a point of each BTZ-line of $\overline N$. 
From Proposition \ref{prop:topology}, $N$ is homeomorphic to $\reg(\overline N)$ and each ends of a Cauchy-surface 
of the latter corresponds to a BTZ-line of $\overline N$. Since  $\mathrm{BTZ-ext}(N)$ contains a point of each BTZ-line of $\overline N$,
 a Cauchy-surface of $\mathrm{BTZ-ext}(N)$ has no ends and is thus compact. 
 From Theorem \ref{theo:BTZ_Cauchy-completeness},  $\mathrm{BTZ-ext}(N)$ is Cauchy-maximal. Thus, from Proposition 
 \ref{prop:c-compact_maximal},  $\mathrm{BTZ-ext}(N)$ is absolutely maximal and $(iii)\Rightarrow (ii)$.

\end{itemize}

\subsubsection{Maximal BTZ-extension of absolutely maximal spacetimes}

Given a regular domain $\Omega$ invariant under the action of some discrete torsionfree subgroup of isometries
$G\subset \isom(\E^{1,2})$.
From the construction above, the only remaining question is whether the line of fixed points of a parabolic element of 
$G$ is in the boundary of $\Omega$.

\begin{prop} \label{prop:fix_line}
 Let $\Omega$ be a regular domain invariant under the action a discrete torsionfree subgroup $G\subset \isom(\E^{1,2})$ such  
 that $\Omega/G$ is absolutely maximal. 
Then, for all $\psi \in G$ parabolic, the BTZ-line associated to $\psi$ is non empty.
\end{prop}
\begin{proof}[Proof of Proposition \ref{prop:fix_line} for type $(i-iii)$ ]

Among the three first cases of Theorem \ref{theo:barbot_mess},  case $(i)$ does not admit parabolic isometry and is thus trivial.
In cases $(ii)$ and $(iii)$, the group $G$ is either generated by a translation in which case the proposition is trivial or  
by a parabolic isometry $\psi$. Since $\Omega/G$ is absolutely maximal, $\Omega=I^+(\fix(\psi))$ and $\Delta_\psi\neq \emptyset$. 
\end{proof}

Only the case $(iv)$ of Theorem \ref{theo:barbot_mess}
remains in which case $L_{|G}$ is injective and $L(G)$ is discrete. \textbf{We now assume 
that $\mathbf{\Omega/G}$  is a type  (iv) spacetime}.

Let $\psi\in G$, let $\Delta=\fix(\psi)$ and let $G_\Delta:=\stab(\Delta)$. We can assume $\Delta$ goes through the origin $O$ of $\E^{1,2}$ and take some $u$ on $\Delta$ above and distinct from $O$.
For $R\in \R$, defines  the planes $$\Pi_R := \{x\in \E^{1,2} ~|~ \langle x|u \rangle = -R \}$$ so that they are the 
the planes parallel to $\Delta^\perp$ and we have 
$$\forall R\in \R, I^-(\Pi_R)=\left\{\left. x\in \E^{1,2} ~\right|~ \langle x|u \rangle > -R \right\}.$$ 
Since $G$ acts on $\Omega$ torsionfree, $\Delta\cap \Omega=\emptyset$, furthermore $\forall x\in \Omega, J^+(x)\subset \Omega$ 
thus $\Omega\subset I^+(\Delta)=I^+(\Pi_0)$.
Let $$R_0 := \max\{R\in \R ~|~ \Pi_R\cap \Omega=\emptyset\} \quad \mathrm{and} \quad \U_{\lambda,R}:=I^+(\lambda u)\cap I^-(\Pi_R)$$
Let $\mathcal C$ be the cone of future lightlike vector from $O$. Up to some translation  of the origin $O$ along $\Delta$, we can assume that $\mathcal C \cap \Pi_{R_0+1}\subset \Omega$.

We want to find some $\lambda>1$ and $R>R_0$ such that $\U_{\lambda, R}$ is disjoint from its translations by $G\setminus G_\Delta$.

\begin{lem}\label{lem:ineq_1}Let $\lambda>1$ and $R=R_0+1$.
  Let $\phi\in G$ such that  $\phi \U_{\lambda,R} \cap \U_{\lambda,R}\neq \emptyset $, then 
  $$ 0\geq \langle \phi_Lu |u \rangle \geq -\frac{R_0+2}{\lambda} $$
 \end{lem}
\begin{proof}
  To begin with, $\langle \phi_L u |u  \rangle$ is non-positive since $\phi_Lu$ and $u$ are both future pointing.

 If $\phi\in G_\Delta$, then $\phi_Lu=\phi u=u$ and thus $\langle  \phi_Lu|u\rangle=0$.  
 If $\phi\notin G_\Delta$.  we have $\phi.(\lambda  u) \in I^-(\Pi_R)$ thus 
 $$\lambda\langle \phi_L|u\rangle +\langle \tau_\phi|u\rangle > -R $$
Furthermore, there exists a unique $v\in \mathcal C$  such that $\phi_Lv= \alpha v$ with $0< v\leq 1$ and $v\in u$. 
Since $L(G)$ is discrete, the vector $v$ is parallel to $u$  if and only if $\phi \in G_\Delta$. Then, for $\phi \notin G_\Delta$, 
$v \in \Pi_R \cap C \subset \Omega $ and  since $G$ stabilises $\Omega$, $\phi v \in \Omega$
and $\phi\langle v |u\rangle \leq -R_0$. Therefore : 
$$-R_0\geq \langle\phi v|u\rangle =   \alpha\langle v|u \rangle+\langle \tau_\phi |u\rangle = -R+\langle \tau_\phi |u\rangle $$
and then : 
$$1=R-R_0\geq \langle \tau_\phi |u\rangle.$$
We thus have : 
$$\langle \phi_Lu|u\rangle \geq \frac{-R-1}{\lambda}= \frac{-R_0-2}{\lambda} $$
\end{proof}

\begin{lem} Let $\lambda>2$, $R=R_0+1$ and let $\phi\in G$.

Then, there exists $\phi_0=1,\phi_1,\cdots \phi_n \in G$ such that 
$$\phi \U_{\lambda,R} \cap \U_{\lambda,R}\neq \emptyset \Rightarrow  
  \phi \in \bigcup_{i=0}^n G_\Delta \phi_i G_\Delta$$
 
\end{lem}
\begin{proof}

Take some $p^* \in \E^{1,2}$ such that $J^+(O)\cap \Pi_R\cap J^-(p^*)$ contains a fundamental domain of the action of $G_\Delta$ on 
$\mathcal C\cap \Pi_R$.
Then, $J^+(O)\cap J^-(p^*)$ contains a fundamental domain of the action of $G_\Delta$ on 
$\mathcal C\cap \Pi_{R'}$ for every $R' \in ]0,R[$.
From Lemma \ref{lem:ineq_1}, for all $k\in \Z$, $$ 0\geq  \langle L(\psi^k\phi) u |u  \rangle >-\frac{R_0+2}{\lambda} $$
and $k\in \Z$ can be chosen such that $ L(\psi^k\phi)u$ is in the fundamental domain of the action of 
$G_\Delta$  on $\mathcal C$ and thus in  $J^+(O)\cap J^-(p^*)$.
Thus, there exists $k\in \Z$ such that $L(\psi^k \phi). p_* \in J^-(p^*)\cap J^+(O)$. 
The diamond $J^+(O)\cap J^-(p^*)$ 
is compact and from \cite{MR918457} $L(G)u$ is discrete, furthermore  $L_{|G}$ 
is injective,  thus there exists only finitely many 
$[\phi']\in G/ G_\Delta$ such that $L(\phi') u \in J^-(p^*)\cap J^+(O)$.
Let $\{\phi_0=1,\phi_1,\cdots,\phi_n\}$ be a set of representative of these $[\phi']$, then :
$$\phi_L \in \bigcup_{i=1}^n G_\Delta L(\phi_i) G_\Delta$$
which yields the results since $L_{|G}$ is injective and $L(G_\Delta)=G_\Delta$.

\end{proof}

\begin{proof}[Proof of Proposition \ref{prop:fix_line}, type (iv)]
Take $$R= R_0+1 \quad \mathrm{and}\quad \lambda = \max_{i\in [1,n]} \left ( \frac{R_0+1} {-\langle L(\phi_i) u |u  \rangle}\right )+1$$
we obtain $\phi \U_{\lambda,R} \cap \U_{\lambda,R} =\emptyset$ for $\phi \in G\setminus G_\Delta$.

Consider the domain $$\Omega':= \Omega\cup \bigcup_{[\phi] \in G/G_\Delta} \phi \U_{\lambda,R}. $$
$\Omega'$ is globally hyperbolic and $G$-invarariant. Let $\widetilde \Sigma$ be $G$-invariant Cauchy-surface of $\Omega$.
Take $p\in I^-(\Pi_R)\cap \Omega$, 
$I^+(p)\cap \Pi_R\subset \Omega$ is the interior of a parabola inside $\Pi_R$ thus every lightlike line in $\Pi_R$ 
intersects $\Omega$. The surface $\widetilde \Sigma$ is a Cauchy-surface of $\Omega$ and thus intersects every lightlike line 
in $\Pi_R$.
Consider $\Omega'/G_\Delta \subset I^+(\Delta)/G_\Delta \simeq \E^{1,2}_0$ and use 
the cylindrical coordinates of $\E^{1,2}_0$. 
Then, the surface $\widetilde \Sigma/G_\Delta$ is a Cauchy-surface of $\Omega/G_\Delta$ and intersects every vertical line $\{r=R,\theta=\theta_0\}$. 
We can use Lemma 53 in \cite{brunswic_btz_ext} to extend $\widetilde \Sigma/G_\Delta \setminus \{r<R\}$ 
to a metrically complete Cauchy-surface of $\Omega'/G_\Delta$ we call $\Sigma$. The lift $\widetilde \Sigma$ of $\Sigma$ is a metrically 
complete Cauchy-surface of $\Omega'$. Then $\Omega'/G$ is Cauchy-complete and we have a natural injective $\E^{1,2}$-morphism 
$\Omega/G \rightarrow \Omega'/G$. Since $\Omega/G$ is absolutely maximal, this map is onto and $\Omega'=\Omega$. 

Finally, $$\fix(\psi) \cap \partial \Omega =\fix(\psi) \cap \partial \Omega' \neq \emptyset.$$

\end{proof}

\begin{cor}\label{cor:compact}
 The maximal BTZ-extension of an absolutely maximal, Cauchy-complete, $\E^{1,2}$-manifold of admissible holonomy is Cauchy-compact.
\end{cor}
\begin{proof}The holonomy of $M$ is admissible, in particular, $M$ is homeomorphic to $\Sigma^* \times \R$ for some compact surface 
$\Sigma$ and finite number of punctures $S$. 
 Moreover, the holonomy of a peripheral loop around a puncture in $S$, is parabolic and, from Proposition   \ref{prop:fix_line},
 admits an open half-line of fixed point in the boundary of $\Omega$. Then, by Theorem \ref{theo:max_BTZ_ext}, the maximal 
 BTZ-extension of $M$ is $\widetilde{\Omega}(G)/G$, which is homeomorphic to $\Sigma\times \R$ and thus Cauchy-compact.
\end{proof}

\subsection{Construction of spacetimes of given admissible holonomy}\label{sec:admissible_space_time}
Our goal is to construct a globally hyperbolic Cauchy-complete spacetime of given admissible holonomy. 

\begin{prop}\label{prop:tau-susp} Let $\Sigma$ be a compact surface of genus $g$ and let $S$ be a finite subset of cardinal $s>0$ such that $2g-2+s >0$.
Let $\Gamma:=\pi_1(\Sigma^*)$.
  
  The map 
  $$\fonction{\d\susp_{\H^2}}{T\teich_{g,s}}{\M_{g,s} (\E^{1,2})}{\rho}{\Omega(\rho)/\rho} $$
   with $\Omega(\rho)$ defined in Proposition \ref{prop:abs_ext_2}, is well defined and inverse to the holonomy map. 
\end{prop}
The main point is to prove that $\Omega(\rho)$ is non empty.
Let $\Sigma$ be a  compact surface, $S$ a finite subset of $\Sigma$, $\Gamma:=\pi_1(\Sigma^*)$ and let 
$\rho : \Gamma \rightarrow \isom(\E^{1,2}) $ a marked admissible representation. 
We assume that the linear part $\rho_L$ of $\rho$ is discrete and faithful and aim to construct a spacetime of holonomy $\rho$.

We use divergent direction fields to define a locally injective continuous $\Gamma$-invariant map  $\H^2\times \R_+^*\rightarrow \E^{1,2}$ where 
$\Gamma$ acts on trivially on $\R_+^*$, via $\rho_L$ on $\H^2$ and via $\rho$ on $\E^{1,2}$. The image of $\H^2\times \{t\}$ will
be $\Gamma$-invariant acausal and metrically complete surfaces in $\E^{1,2}$. This procedure is taken from a still unpublished 
work of Barbot and Meusburger \cite{spinexp} on construction of spacetimes with particles with spind.

\begin{defi}[Divergent direction field] Let $X\subset \H^2$ be any subset of $\H^2$. 
\begin{itemize}
 \item A direction field on $X$ is a map $ f : X \rightarrow \E^{1,2}$. 
 \item A direction field on $X$ is divergent if for all $x,y\in \H^2$ $$\langle  f(x)-f(y)  | x+\langle x|y\rangle y\rangle\geq 0  $$
 \item A direction field on $X$ is locally divergent  if for all $p\in \H^2$, there exists an open neighborhood $\U$ around $p$ such
 that $f_{|\U}$ is divergent.
\end{itemize}

\end{defi}

\begin{defi} A direction field $f$ on $\H^2$ is $\rho$-equivariant if for all $p\in \H^2$ and all $\gamma \in \Gamma$,
$$f(\rho_L(\gamma)p)=\rho(\gamma)f(p).$$
 
\end{defi}

\begin{defi} Let $f$ be a direction field.
  Define the map $$\fonction{\D_f}{\H^2 \times \R_+^*}{ \E^{1,2}}{(x,t)}{f(x)+t x}$$
\end{defi}

\begin{rem} Let $\rho : \Gamma\rightarrow \isom(\E^{1,2})$ a morphism.
 If $f$ is a locally divergent and $\rho$-equivariant direction field then there exists a unique $\E^{1,2}$-structure on $\H^2/\Gamma \times \R_+^*$ such that 
 $\D_f$ is its developping map. 
 Furthermore, the holonomy associated to $\D_f$ is $\rho$.
\end{rem}

\begin{lem} \label{lem:div_C1} Let $X\subset \H^2$ be an open subset and let $f:X\rightarrow \E^{1,2}$ be a $\C^1$ direction field on $X$.
If for all $\xi\in T\H^2$, $\langle \d f(\xi)| \xi\rangle > 0$
then $f$ is locally divergent.
\end{lem}
\begin{proof} 
Restricting to smaller open subsets, we can assume $X$ convex.
 Let $u,v$ be two element of $X$, and let $\xi$ the unique element of $T_uX$ whose image by the exponential is $v$. We have :
 $$v=\frac{\sinh(\|\xi\|)}{\|\xi\|}\xi +\cosh(\|\xi\|)u $$
 and, since $\langle u |\xi \rangle = 0$ : 
 \begin{eqnarray}
  v+\langle u |v\rangle &=& \frac{\sinh(\|\xi\|)}{\|\xi\|}\xi +\cosh(\|\xi\|)-\cosh(\|\xi\|)u \\
  &=&  \frac{\sinh(\|\xi\|)}{\|\xi\|}\xi 
 \end{eqnarray}
We also have: 
\begin{eqnarray}
 \|u-v\|^2&=& \sinh(\|\xi\|)^2-\left(\cosh(\|\xi\|)-1\right)^2 \\
 &=& 2\cosh(\|\xi\|)-2\\
 &=& 4\sinh^2(\|\xi\|/2)
\end{eqnarray}

Since $f$ is $\C^1$, there is a continuous map $\varepsilon : X\times X \rightarrow \E^{1,2}$, vanishing on the diagonal, such that:
$$f(v)=f(u)+\d f(\xi)+\|v-u\| \varepsilon(u,v) $$
Then: 
\begin{eqnarray}
 \langle f(v)-f(u) | \langle u|v\rangle u\rangle &=& \left\langle f(u)+\d f(\xi)+\|v-u\| \varepsilon(u,v) \left|\frac{\sinh(\|\xi\|)}{\|\xi\|}\xi \right. \right\rangle \\
 &=& \frac{\sinh(\|\xi\|)}{\|\xi\|} \langle \d f(\xi)|\xi\rangle + 2\sinh(\|\xi\|/2)\langle \varepsilon(u,v)|\xi\rangle
\end{eqnarray}

Let $u_0\in X$, let $\varepsilon_0 = \frac{1}{2} \min_{\xi \in T_{u_0}X,\|xi\|=1} \langle \d f(\xi)|\xi\rangle$.
Since $\varepsilon$ and $\d f$ are positive and continuous on $X\times X$, there exists a 
neighborhood $U$ of $u_0$ such that for all $(u,v)\in \U^2$, $\langle \varepsilon(u,v)|\xi\rangle \leq \frac{\varepsilon_0}{2} \|\xi\|$ and $\langle\d f(\xi)|\xi\rangle\geq \|\xi\|^2\varepsilon_0$.
Then for all $(u,v)\in\U^2$, $\langle f(v)-f(u) | \langle u|v\rangle u\rangle \geq0$, thus $f$ is divergent on $\U$.
\end{proof}

\begin{prop}\label{prop:direction_field} There exists a $\rho$-equivariant locally divergent field $g$ on $\H^2$ such that 
 for all $t\in \R_+^*$, $\D_g(\H^2,t)$ is acausal and  complete.

\end{prop}

\begin{proof}
Consider $(\gamma T_{i})_{\gamma \in \Gamma,i\in [1,n]}$ the lift of an ideal geodesic triangulation of $\H^2/\Gamma$. 
Take a disjoint family of embedded horodisks $(\gamma H_j)_{j\in\{1,\cdots,s\},\gamma\in \Gamma}$, 
we get a non-geodesic cellulation $(\gamma T_i',H_j)_{i\in [1,n],j\in [1,s]}$ where $T_i':= T_i\setminus \bigcup_{j=1}^s H_j$.

Since $\rho$ is admissible, for every $j\in [1,s]$, there exists a unique 
$\Delta_j $  such that all $\gamma$ parabolic fixing $H_j$ set-wise fixes $\Delta_j$ point-wise.
For each $j$, choose a point $p_j\in \Delta_j$ and set 
$$\forall x\in \gamma H_j, f(x)=\rho(\gamma) p_j.$$ 
Let $\varphi:[0,1]\rightarrow [0,1]$, $\C^{1}$, inscreasing and such that $\varphi(0)=\varphi'(0)=\varphi'(1)=0$
and $\varphi(1)=1$ and For $i\in I$, the cell $T_i'$ is a non-geodesic hexagon $[A_k]_{k\in \Z/6\Z}$ with $\C^1_{pw}$ boundary. 
 $f$ is defined on the horocycle part $[A_0A_1]\cup [A_2A_3]\cup[A_4A_5]$.
Extend $f$ on $[A_{2k+1}A_{2k+2}], k\in\{0,1,2\}$ putting 
$$f(x)=\varphi\left(\frac{d_{\H^2}(x,A_{2k+1})}{d_{\H^2}(A_{2k+1},A_{2k+2})}\right) f(A_{2k+2}) + (1-\varphi)\left(\frac{d_{\H^2}(x,A_{2k+1})}{d_{\H^2}(A_{2k+1},A_{2k+2})}\right)f(A_{2k+1})$$
Then extend $f$ in a $\C^1$ way to $T_i'$ in such a way that $\d_x f.h=0$ for  $x\in\partial T_i'$ and $h\perp \partial T_i'$.
This way, $f$ is a $\C^1$, $\Gamma$-invariant direction field on $\H^2$. Then, $\xi \mapsto, 
\langle\d f(\gamma \xi)|\gamma \xi\rangle=\langle\d f(\xi)| \xi\rangle$, is $\Gamma$-invariant (with $\Gamma$ acting trivially on $\R$),
homogeneous of degree 2 on each fiber and zero in every $TH_j$, $j\in[1,s]$. Therefore, writing $T^1(\H^2/\Gamma)$ the unitary fiber 
bundle over $\H^2/\Gamma$,
$$\fonctionn{T^1(\H^2/\Gamma)}{\R}{\xi}{\langle\d f(\xi)| \xi\rangle} $$
is well defined, have compact support and is thus bounded, let $M\in \R$ be its minimum. 
The function $g:x\mapsto f(x)+(M+1)x$ then satisfies the hypothesis of Lemma \ref{lem:div_C1} and is thus locally divergent.

Let  $t \in \R_+^*$, the quadratic form induced on the level set  $\D_g(\H^2,t)$ by metric of $\E^{1,2}$
is $q(\xi)=\langle \d (g+t \mathrm{Id})(\xi)|\xi\rangle=\langle \d g(\xi)|\xi\rangle+t\|\xi\|^2\geq (1+t)\|\xi\|^2>0$. Then, this level set 
is a closed spacelike hypersurface of $\E^{1,2}$. By Corollary 46, chapter 14 of \cite{oneil}, $\D_g(\H^2,t)$ is acausal. 
Since $q(\xi)\geq (1+t)\|\xi\|$, $x\mapsto \D_g(x,t)$ enlarges distances, and since $\H^2$ is metrically complete, $\D_g(\H^2,t)$ is metrically complete. 
\end{proof}

\begin{proof}[Proof of Proposition \ref{prop:tau-susp}]
 Let $g$ given by Proposition \ref{prop:direction_field}, let $\Omega$ be the Cauchy development of $\D_g(\H^2,1)$ in $\E^{1,2}$ and let 
 $M:= \Omega/\Gamma$. 
 The spacetime $M$ is globally hyperbolic, Cauchy-complete, Cauchy-maximal, future complete, homeomorphic to $\H^2/\Gamma \times \R$ 
 and its holonomy is $\rho$. By Proposition \ref{prop:abs_ext_2}, $\Omega(\rho)$ is not empty and 
 $\Omega(\rho)/\rho$ is the absolutely maximal extension of $M$ .
\end{proof}

\begin{rem} The introduction of divergent direction field could have been avoided since we only needed a $\Gamma$-invariant 
metrically complete and acausal surface in Minkowski space. However, the procedure is way more general and allows to construct 
globally hyperbolic spacetimes of holonomy $\rho$ which is may not be admissible. 
 
\end{rem}

\subsection{$\H^2-\E^{1,2}_0$ correspondances}\label{sec:admissible_hol}

The last section ended with the definition of the map $\d\susp_{\H^2}$. We can thus state and prove the following Theorem.
\begin{theo} \label{theo:Teich}
The following maps are well defined and bijective.

$$
\xymatrix{
\teich_{g,s} \ar@<1ex>[rr]^{\susp_{\H^2}}&&  \ar@<1ex>[ll]^{\susp^{-1}_ {\H^2}} \M_{g,s}^L(\E^{1,2)}  \ar@<1ex>[rr]^{\mathrm{BTZ-ext}}&&  \ar@<1ex>[ll]^{\mathrm{Reg}}\M_{g,s}(\E^{1,2}_0) 
}
$$ $$
\xymatrix{
T\teich_{g,s} \ar@<1ex>[rr]^{\d\susp_{\H^2}}&&  \ar@<1ex>[ll]^{\Hol}\M_{g,s}(\E^{1,2})  \ar@<1ex>[rr]^{\mathrm{BTZ-ext}}&&  \ar@<1ex>[ll]^{\mathrm{Reg}} \M_{g,s}(\E^{1,2}_0) 
}
$$
\end{theo}

\begin{proof}
 Proposition \ref{prop:tau-susp} shows that the map $\d\susp_{\H^2}$ is bijective and the  inverse of $\Hol$. 
 Corollary \ref{cor:compact} shows that the maximal BTZ-extension of an absolutely maximal Cauchy-complete $\E^{1,2}$-manifold 
 of admissible holonomy is Cauchy-compact. Thus the map $\mathrm{BTZ-ext}$ is well defined. 
 
 Consider $M$ a Cauchy-compact globally hyperbolic $\E^{1,2}_0$-manifold, from Theorem 
 \ref{theo:BTZ_Cauchy-completeness} $M':=\reg(M)$ is Cauchy-maximal and Cauchy-complete.
 By Lemma \ref{lem:cauchy_max_abs_max}, 
 $M$ is absolutely maximal and by Lemma \ref{lem:abs_max_compact}  $M'$ is absolutely 
 maximal. 
 Assume $M = \Sigma\times \R$ with $\sing_0(M)=S\times \R$, note $\rho:\Gamma \rightarrow \isom(\E^{1,2})$ its holonomy and $\Omega$ the
 developpement of $M'$. The group $\Gamma$ acts properly discontinuously and freely on $\Omega$ 
 thus $\rho(\Gamma)$ has no elliptic element.
 Consider a peripheral loop $\gamma$ around some puncture in $S$, it can be chosen in a neighborhood chart of a BTZ point 
 and thus $\rho(\gamma)$ is parabolic. 
 Let $\gamma \in \Gamma$ such that $\rho(\gamma)$ is parabolic, then, as proved in Corollary 
 \ref{cor:admissible},  $\rho(\gamma)$ has a line of fixed
 point $\Delta$. Since  $M'$ is  absolutely maxima Lemma \ref{prop:fix_line} shows that
 $\Delta$  is adjacent to $\Omega$. Then, $\Delta$ clorresponds to a BTZ-line in the maximal BTZ-extension of $\reg(M)$ 
 which  is $M$ and $\rho(\gamma)\in\stab(\Delta)=\rho(\langle c \rangle )$ for some $c$ peripheral.
 Since $\rho$ is faithful, $\gamma\in \langle c\rangle$ thus $\gamma$ is peripheral. 
 Finally, $\rho_L$ is admissible and from Corollary \ref{cor:admissible} so is $\rho$. 

  The map $\reg$ is then well defined and since $\mathrm{BTZ-ext}$ and $\reg$ are inverse of each other, both are bijective.
 
\end{proof}

\section{Decorated Moduli correspondances }
\subsection{Decorated Moduli spaces and Penner surface}\label{sec:Penner}
Let $g \in \N$, $s\in \N$ and let $\Sigma$ be a compact surface of genus $g$. In this section, $(\Sigma,S)$ denote a topological
surface $\Sigma$ together with a set of $s$ marked point $S$.

In \cite{penner1987}, Penner introduced a so-called Decorated Teichmüller space $\widetilde \teich_{g,s}$. Penner defines it as the fiber bundle 
over the Teichmüller space which points are  equivalence classes of marked hyperbolic surface of finite volume  with a choice of an horocycle around each cusp.
Penner then construct a polyhedral surface associated to a point of his decorated Teichmüller space. 
The construction described in section 1-4 of \cite{penner1987} goes as follows. 

Let $\Sigma^*$  be a hyperbolic surface of genus $g$ with $s$ cusps and choose an horocycle $h_i$ around each cusp. 
The universal cover of $\Sigma^*$ is identified with the hyperbolic plane embedded into Minkowski space. Notice that 
a  future light ray  from the origin $O$ of Minkoswki corresponds to a point at infinity in the boundary of $\H^2$.
Then an horoball on $\H^2$ centered on a point $r\in \partial \H^2$ is exactly the intersection $J^+(p)\cap \H^2$ for some $p$ on the light ray 
corresponding to $r$.
Thus to each $h_i$ corresponds a unique point $p$ on the future light cone from the origin. 
The idea is then to take the boundary of the closed convex hull of the points corresponding to the horocycles $h_i$.
Penner prove that this boundary, say $\Sigma'$, is polyhedral in the sens that there is locally a finite number of totally geodesic
2-facets around each point of $I^+(O)$. He also proves the 2-facets are all spacelike 
and that each future timelike ray from the origin intersects $\Sigma'$ exactly once.  Moreover, the vertices are 
all in the future light cone $\partial J^+(O)$. 
Since $\Sigma'$ is $\pi_1(\Sigma^*)$-invariant, this means that $\Sigma'/\pi_1(\Sigma^*)$ is naturally endowed with a singular euclidean metric.

In our framework, the first step is exactly the suspension and the second consists in taking the maximal BTZ-extension 
and then associate to each horocycle a point of the BTZ-lines.
We thus define a Lorentzian analogue of Penner decorated Teichmüller space by defining decoration of a $\E^{1,2}_0$-manifold.

\begin{defi}[Decoration of $\E^{1,2}_0$-manifold] Let $M$ be $\E^{1,2}_0$-manifold. A decoration of $M$ 
is a choice of a point $p_{\Delta}$ on each connected BTZ line $\Delta$ of $M$.
\end{defi}

\begin{defi}
Let $(M_1,h_1,p^1_{1},\cdots,p^1_{s} ), (M_2,h_2,p^2_{1},\cdots,p^2_{s})$  be two decorated marked $\E^{1,2}_0$-manifolds homeomorphic 
to $\Sigma\times \R$ and such that for $i=1,2$  $\sing_0(M_i)$ has exactly $s$ connected components
decorated by the points $(p_j^i)_{j\in[1,s]}$. 

They are equivalent if there exists a $\E^{1,2}_0$-isomorphism $\varphi:M_1\rightarrow M_2$ such that $\varphi(p_i^1)= p_i^{2}$ for $i\in [1,s]$ and such 
that $h_2^{-1}\circ  \varphi\circ h_1 $ is homotopic to $Id_{\Sigma\times \R}$.
\end{defi}

\begin{defi}[Decorated BTZ Moduli space]
Let $g\in \N$, let $s\in \N$ and let $\Sigma$ be a compact surface of genus $g$.

 The decorated $\E^{1,2}_0$-moduli space $\widetilde \M_{g,s}(\E^{1,2}_0)$ is 
 the  set of equivalence classes of decorated marked
$\E^{1,2}_0$-manifold homeomorphic to $\Sigma\times \R$ and such that $\sing_0(M)$ has exactly $s$ connected components.  
\end{defi}

\begin{rem}
We can define accordingly the Decorated linear BTZ Moduli space $\widetilde \M_{g,s}^L(\E^{1,2}_0)$. 
\end{rem}
\begin{defi} 
Let $\mathcal H$ be the set of horocycle of $\H^2\subset \E^{1,2}$.
$$\fonction{\dec^{-1}}{\partial J^+(O)}{\mathcal H}{p}{\partial J^+(p)\cap \H^2} $$ 
\end{defi}
\begin{lem}[\cite{penner1987}]\label{lem:dec}
  The map $\dec^{-1}$ is bijective and continuous.
\end{lem}

Write $\dec$ the inverse  of $\dec^{-1}$.

\begin{prop}
The following maps are bijective :
 $$\xymatrix{ \widetilde \teich_{g,s} \ar@<1ex>[rrr]^{(\mathrm{BTZ-ext} ~\circ~ \susp_{\H^2})\oplus \dec}&&&\ar@<1ex>[lll]^{(\susp_{\H^2}^{-1}~\circ~\reg)\oplus \dec^{-1}} \widetilde \M_{g,s}^L(\E^{1,2}_0)} $$  
$$\xymatrix{ \widetilde {T\teich}_{g,s} \ar@<1ex>[rrr]^{(\mathrm{BTZ-ext} ~\circ~ \d\susp_{\H^2})\oplus \dec}&
&&\ar@<1ex>[lll]^{\Hol~\oplus~\dec^{-1}} \widetilde \M_{g,s}(\E^{1,2}_0)} $$  
 \end{prop}

\begin{proof}
 It follows from Theorem \ref{theo:Teich} and \ref{lem:dec}.
\end{proof}

Equivalence between marked singular euclidean surfaces is defined is a similar manner as for the other equivalences introduced so far.
Then, we can define the Euclidean moduli space.

\begin{defi} 
The Euclidean moduli space $\M_{g,s}(\E^{2})$ is the space of equivalence classes of marked singular euclidean 
surface homeomorphic to $\Sigma$ with exactly $s$ conical singularities.
\end{defi}

The second part of the construction of Penner then defines a map 
$$\xymatrix{ \widetilde \M_{g,s}^L(\E^{1,2}_0)\ar[rr]^{\pen}&&\M_{g,s}(\E^{2}) }$$
Our objectives now are to construct an inverse to  the map $\pen$ and to extend $\pen$ to $\widetilde \M_{g,s}(\E^{1,2}_0)$.

\subsection{Suspension of a singular euclidean surface}

The first step of our construction of the inverse of Penner map is to describe the cellulation which comes 
with the Penner-Epstein surface.
First, 
Proposition 2.6 of \cite{penner1987} shows  
that the developpement of a cell of the Penner surface satisfies Ptoleme equality and thus its vertices are cocyclic.
Second, it shows the a Ptoleme inequality holds for every quadrilateral about an edge of the cellulation.
This is a caracterisation of the  Delaunay cellulation of a compact  singular Euclidean surface we now describe.

\begin{defi} A geodesic cellulation of $(\Sigma,S)$ is ideal if the edges start from $S$ and ends on $S$.  
\end{defi}

\begin{defi}[Hinge]  A hinge is a euclidean quadrilateral together with a diagonal.
\end{defi}
\begin{defi}[Hinge about an edge]
Let $\mathcal T$ be an ideal geodesic triangulation of $(\Sigma,S)$. 

The hinge about an edge $e$ of the triangulation $\mathcal T$ is the hinge
obtained by gluing along $e$ the two triangles of $\T$ intersecting along $e$.
The diagonal of this hinge is the edge $e$. 
\end{defi}
\begin{defi}[Legal edge]
Let $\mathcal T$ be an ideal geodesic triangulation of $(\Sigma,S)$

An edge $e$ of $\mathcal T$ is legal if the hinge $h$ about $e$ is not inside the circumscribed circle around one of the triangles composing $h$.
\end{defi}

\begin{prop}[\cite{Troyanov}] Let $(\Sigma,S)$ be a singular euclidean surface.  There exists a unique ideal 
 celluation $D$ such that every edge of $D$ is legal. Moreover, every cell of $D$ is a cocyclic polygon.
\end{prop}

The second step of our construction of the inverse of Penner map is to use Lemma 2.4 and Corollary 2.5 of \cite{penner1987}. we rewrite them slightly 
to convey our needs. 
\begin{lem}[\cite{penner1987}]\label{lem:inscription} Let $C$ be a euclidean cocyclic polygon. Then there exists a totally geodesic embedding of $C$ into $\E^{1,2}$ such that 
the vertices of $C$ are in the future light cone $\partial J^+(O)$. 

Moreover, two such totally geodesic embeddings only differ by an isometry $\gamma\in \SO_0(1,2)$. 
\end{lem}

\begin{defi}[Suspension of a cocyclic Euclidean polygon]
Let $C$ a cocyclic Euclidean polygon and let $i: C \rightarrow \E^{1,2}$ be the a totally geodesic embedding of $C$
such that the vertices of $i(C)$ are in the light cone $\partial J^+(O)$. Write $q(t,x,y)=-t^2+x^2+y^2$.

  The suspension of a cocyclic Euclidean polygon $C$ is 
	$$ \susp_{\E^2}(C):= (C\times \R_+^*,\d s^2) \quad \mathrm{with}\quad   \d s^2 (x)= q(i(x))\d t^2 + \d s^2_{C}(x)$$
	together with the decoration $\left(i(v_i)\right)_{i\in [1,k]}$ where $v_1,\cdots,v_k$ are the vertices of $C$.
\end{defi}

\begin{rem} The suspension of $C$ is nothing more than the metric on the cone or rays from $O$ intersecting a totally geodesic
embedding of $C$ enscribed in the future light cone from the origin.
\end{rem}
\begin{rem} This suspension does not depend on the choice of the embedding since from Lemma \ref{lem:inscription}  all of them are isometric via a global isometry 
of $\E^{1,2}$ .
\end{rem}

The next step in our construction of the inverse of Penner map is to glue the suspension of each cell of the Delaunay 
cellulation of the singular Euclidean surface $(\Sigma,S)$. In the same manner as for gluing of Euclidean triangles, 
gluings give rise to singularities. The next section is devoted to the proof that these singularities are BTZ lines. 
If the reader is convinced of this fact, he may skip the following section.

\subsection{Gluings of future lightlike Minkowski wedge} \label{sec:BTZ_gluing}

This section details general properties of gluings of Minkowski wedges (see Definition \ref{def:wedge} below). 
The following proofs are classical in the riemannian contex and are not much more complicated in the Lorentzian context.
However, the list of singularties that arise in the Lorentzian context is
longer than in the Riemannian context,  
a complete classification with more involved properties is describred in \cite{Particles_1,Particles_2}. 
Since we are only interested in gluings of direct future wedges, the aim of the section is Proposition \ref{prop:gluing_BTZ}
which shows such gluings only give rise to BTZ singularities. 

\begin{defi}[Future wedge in Minkowski] \label{def:wedge}
Let $\Delta$ be a causal line in $\E^{1,2}$ and let $\Pi_1$ and $\Pi_2$ be two non parallel half (1,1)-planes in $\E^{1,2}$ such that $\partial \Pi_1  =\partial \Pi_1 =\Delta$ .
Assume $\Pi_i \in J^+(\Delta)$, then the convex hull of $\Pi_1 \cup \Pi_2$ is a futur wedge in Minkowski space of axis $\Delta$. 

A wedge is oriented by taking a futur vector $u_0$ in $\Delta$ and a vectors $u_i \in \Pi_i, i \in \{1,2\}$. Its orientation is direct if the basis 
$(u_0,u_1,u_2)$ is a direct in $\E^{1,2}$. A wedge directly oriented is said direct.
\end{defi}

\begin{lem}\label{lem:unique_struc_0}
 For every future direct wedges $S$ and $S'$, there exists $\gamma \in \isom(\E^{1,2})$ such that $S'=\gamma S$
\end{lem}
\begin{proof}
 Applying some translation, we can choose the origin of $\E^{1,2}$ to be on $\Delta$. 
 Then $\Delta$ is a point on $\partial \H^2$ and the intersection of $\Pi_1$ and 
 $\Pi_2$ with $\H_2$ are geodesics intersecting on the boundary at $\Delta$. Together with the other intersection point of $\Pi_1$  and 
 $\Pi_2$ with $\partial \H^2$ we get a triplet of points $(\Delta,A_1,A_2)$ on $\H^2$ which totally characterises $S$. 
 The direct condition is equivalent to the fact that $\Delta<A_1<A_2<\Delta $ for the direct orientation on $\partial \H^2$.
 The action of $\SO_0(1,2)$ is 3-transitive on $\partial \H^2$ thus, given two wedges $(\Delta,A_1,A_2)$ and $(\Delta',A_1',A_2')$, 
 then there exists $\gamma$ sending the one to the other.
 
\end{proof}

\begin{defi}[Direct wedge Gluing]
  Let  $S=(\Pi_{1},\Pi_{2})$  and $S'=(\Pi_{1},\Pi_{2})$ be two direct future wedges of respective axis $\Delta$ and $\Delta'$. 
  A direct identification of $S$ to $S'$ is an isometry $\gamma \in \E^{1,2}$ sending $\Pi_{2}$ on  $\Pi_{1}'$ and $\Delta$ on $\Delta'$.
  The gluing of $S$ and $S'$ along $\gamma$ is $$S\oplus^\gamma S' := (S\cup S')/\gamma$$
  
  \end{defi}
  \begin{defi}[Direct cyclic wedge gluing]
  Let $n\in \N^*$, let $S^{(i)}=(\Pi_{1}^{(i)},\Pi_{2}^{(i)})$, for $i\in\Z/n\Z$, be a family of direct future wedges 
  of respective axis $\Delta^{(i)}$ and let $\gamma\in\isom(\E^{1,2})^n$ such that for $i\in [1,n]$,
  $\gamma_i \Pi_2^{(i)} = \Pi_1^{(i+1)}$. The direct gluing of $(S,\gamma)$ is 
  $$\bigoplus_{i\in \Z/\Z}^\gamma S^{(i)} := \left(\bigcup_{i\in \Z/n\Z} S^{(i)}\right)/\sim \quad \quad \mathrm{with}  \quad\quad x\sim y ~\Leftrightarrow~ \exists i\in \Z/n\Z, \gamma_i x= y ~\mathrm{or}~ \gamma_i y= x   $$
\end{defi}

\begin{lem} \label{lem:unique_struc_1}
 There exists a unique $\E^{1,2}$-structure on $(S\oplus_\gamma S')\setminus \Delta$ which extends the $\E^{1,2}$-structure on $Int(S) \cup Int(S')$.
 
 Moreover, $S\oplus_\gamma S'$ is isomorphic to a direct future wedge.
\end{lem}
\begin{proof}
To begin with, the identification and the wedges are direct,
then the quotient only identifies point of $\Pi_2$ to points on $\Pi_1'$.
The natural projection $\pi:S\cup S'\rightarrow S\oplus_\gamma S'$  restricted to the interior of $S$ and $S'$ 
is an homeomorphism onto its image and thus there is a natural $\E^{1,2}$-structure on the image of the interior of $S$ and $S'$. 
Notice that $S\oplus_\gamma S'$ is simply connected.

\begin{itemize}
\item Define  : 
  $$\fonction{\D}{S\oplus_\gamma S'}{\E^{1,2}}{\pi(x)}{\left\{\begin{array}{ll}
                                            x & \mathrm{if}~x\in S'\\ \gamma x& \mathrm{if}~x\in S
                                           \end{array}\right.
} $$
Since the gluing is direct, $\D$ is injective and thus is a local homeomorphism. The pull back of the $\E^{1,2}$-structure 
of $\E^{1,2}$ via $\D$ defines a $\E^{1,2}$-structure on $S\oplus_\gamma S'$ extending the one on the interior of $S$ and $S'$.
The image of $\D$ is exactly the future direct wedge $(\Pi_1,\Pi_2')$

 \item Assume there exists a $\E^{1,2}$-structure on $S\oplus_\gamma S'$ extending the one on the interior of the wedges. 
Let $\D'$ be a developpement of $S\oplus_\gamma S'$ endowed with such an $\E^{1,2}$.
Then the developpement of the interior of $S$ is $\alpha S$ for some $\alpha\in\isom(\E^{1,2})$ and the developpement 
of $S'$ is $\alpha' S$ for some $\alpha'\in \isom(\E^{1,2})$. We choose $\D$ such that $\alpha'=1$.

Considering a curve $c:]-1,1[\rightarrow S\oplus_\gamma S'$ such 
that $c(0) =\pi(x)$ for some $x\in \Pi_2$ then we see that 
$$\forall x\in \Pi_2,\quad  \alpha x = \D'(c(0))=\lim_{t\rightarrow 0^-} \D'( c(t))= \lim_{t\rightarrow 0^+ } \D'(c(t))=\D'(c(0))=\alpha' \gamma =\gamma x$$
thus $\alpha=\gamma$. This proves that  $\D'=\D$.
\end{itemize}
\end{proof}
\begin{lem} \label{lem:unique_struc_2} Let $n\in \N^*$, let $\left(S^{(i)}\right)_{i\in \Z/n\Z}$ be a family of direct future wedges and let $\gamma\in \isom(\E^{1,2})^n$ 
be an identification. 

Then there exists a unique $\E^{1,2}$-structure on the cyclic gluing of $S$ along $\gamma$ extending 
the $\E^{1,2}$-structure on the interior of the $S^{(i)}$.
\end{lem}
\begin{proof}
 Proceed the same as for Lemma  \ref{lem:unique_struc_1}.
\end{proof}

\begin{prop}\label{prop:gluing_BTZ}
 Let $n\in \N^*$, let $(S^{(i)})_{i\in \Z/n\Z}$ be a finite family of direct futur wedges in $\E^{1,2}$ and let 
 $\gamma \in \isom(\E^{1,2})^n$ an identification.

 Then there exists a unique $\E^{1,2}_0$-structure on the cyclic gluing of $S$ via $\gamma$ 
 extending the $\E^{1,2}$-structure of  the interior of $S^{(i)}, i\in\Z/n\Z$.
 Furthermore :
 \begin{enumerate}
  \item $\sing_0(M)=\Delta$
  \item $\Delta$ is a BTZ-line
 \end{enumerate}
\end{prop}

\begin{proof}
 From Lemmas \ref{lem:unique_struc_0}, \ref{lem:unique_struc_1} and \ref{lem:unique_struc_2}, the problem reduces to $n=1$, $S=(\Pi_1,\Pi_2)$ of axis $\Delta$
 and $\gamma$ such that $\gamma \Pi_1 = \Pi_2$ and $\gamma$ fixes point-wise $\Delta$. The isometry $\gamma$ is then 
 parabolic. Let $M$ be the cyclic gluing of $S$ via $\gamma$, then let $\widetilde M := \bigoplus^{(\gamma^{n})_{n\in \Z}}_{n\in\Z} \gamma^n S $ 
 be the non-cyclic gluings of the iterated of $S$ via $\gamma$. Let $\Z$ acts via left multiplication by $\gamma$ on 
 $\widetilde M$.
 The identity on $\bigcup_{n\in \Z} \gamma^n S \rightarrow \E^{1,2}$ quotient out to give a 
 $\Z$-invariant map $\pi : \widetilde M \rightarrow M$ and bijection 
 $\widetilde \D : \widetilde M \rightarrow J^+(\Delta)$.  
 Be careful that $\pi$ is not an homeomorphism, but we can check that the pull back via $\pi\circ \D^{-1}$ of the topology of $M$
 is exactly the BTZ-topology on $J^+(\Delta)$. 
 Thus by quotienting out, we obtain an homeomorphism 
 $$\xymatrix{J^+(\Delta)/\langle\gamma\rangle\simeq \E^{1,2}_0 \ar[rrr]^{~~~~~~~\overline{\pi\circ \D^{-1}} }&&& M}$$
 
 This homeomorphism is a $\E^{1,2}$-morphism on the interior of $\D(S) \subset I^+(\Delta)$, by Lemma \ref{lem:unique_struc_2}, 
 the image of the $\E^{1,2}$-structure on $\reg(\E^{1,2}_0)$ is thus the natural $\E^{1,2}$-structure of $M\setminus \Delta$.
 Then for any $\E^{1,2}_0$-structure on $M$, $\overline{\pi\circ \D^{-1}} $ is a $\E^{1,2}_0$-isomorphism 
 and we can define one as the $\E^{1,2}_0$-structure image from $\overline{\pi\circ \D^{-1}} $.
 
\end{proof}

\subsection{Inverse of Penner map}
We now conclude our construction of the Penner map. 
\begin{lem}
Let $C_1$ and $C_2$ be two cones in $\E^{1,2}$ from $O$ of respective direct triangular basis $[A_1B_1C_1]$ and $[A_2B_2C_2]$. 

Then 
 there exists a unique isometry $\gamma \in \SO_0(1,2)$ such that $\gamma A_1=B_2$, $\gamma B_1=A_2 $ and such that 
 $\gamma C_1$ and 
 $C_2$ are on different sides of the plane $(OB_2A_2)$.
\end{lem}
\begin{proof}
 It a direct Corollary of Lemma 2.3 in \cite{penner1987}.
\end{proof}

\begin{cor} Let $C=[v_1\cdots v_p]$ and $C'=[v'_1\cdots v'_q]$ two direct Euclidean cocyclic 
polygons such that lengths $v_1v_2$ and $v_1'v_2'$ are equal. 
 
 Then, there exists a unique $\gamma \in \isom(\E^{1,2})$ sending the vertex of
 $\susp_{\E^2}(C)$ on the vertex of $\susp_{\E^2}(C')$  and such that $\gamma v_1=v'_2$, $\gamma v_2 = v_1'$
 and such that $C$ and $C'$ lie on different sides of the plane $(Ov_2'v_1')$.
\end{cor}

\begin{defi}[Suspension of a Euclidean surface]
 Let $(\Sigma,S)$ be Euclidean surface singular exactly on $S$, let $C^{(i)}=[v_j^{(i)} : j=1, \cdots, p_i]$ be the Delaunay cells
 of $(\Sigma,S)$. 
 
 Define $\susp_{\E^2}(\Sigma)$, the suspension of $\Sigma$, 
 as the unique gluing of the suspension of the Delaunay cells $C^{(i)}$ 
 which sends the decoration $v_a^{(b)}$ on the marking $v_{c}^{(d)}$ whenever $v_a^{(b)}$ and $v_c^{(d)}$ are equal in $\Sigma$.
 It comes with the decoration induced by the decoration of the suspension of each cell.
\end{defi}

\begin{theo} \label{theo:teichlinear}The following maps are bijective.
  $$\xymatrix{ \widetilde \teich_{g,s}\ar@<1ex>[rrrr]^{(\mathrm{BTZ-ext}~\circ~ \susp_{\H^2})\oplus \dec} &&&&\ar@<1ex>[llll]^{( \susp_{\H^2}^{-1}~\circ~ \reg )\oplus \dec^{-1}} \widetilde \M_{g,s}^L(\E^{1,2})\ar@<1ex>[rr]^{\pen}&&\ar@<1ex>[ll]^{\susp_{\E^2}} \M_{g,s}(\E^2)}$$
\end{theo}
\begin{proof}
 
 By the definition of the maps $\susp_{\E^2}$ and $\pen$ and since the cellulation of the Penner surface is the Delaunay cellulation, 
 $\susp_{\E^2} \circ \pen $ is the identity on $ \teich_{g,s}$.

 Let $[(\Sigma,S)]$ be a point of $\M_{g,s}(\E^2)$.
 By construction, $M:=\susp_{\E^2}(\Sigma)$ is a point of $\widetilde \M_{g,s}^L(\E^{1,2})$.
 The universal cover of $\reg(\susp_{\E^2}(\Sigma))$ is then $I^+(O)$ and the associated augmented regular domain 
 is $I^+(O)$ together with a set of lightlike rays $\Delta$ from $O$. Each of the ray corresponds to a BTZ line of $\susp_{\E^2}(\Sigma)$
 and the decoration gives a family of points $p_\Delta$ on each $\Delta$. The totally geodesic embedding 
 of $\widetilde{\Sigma^*}$ in $I^+(O)$ given by the suspension is a polyhedral surface. Furthermore, 
 every edge of the Delaunay cellulation of $\Sigma$ is legal and thus the hinge about an edge 
 satisfies the Ptoleme inequality. Then, from Lemma 2.6 of \cite{penner1987}, the totally geodesic embedding of 
 $\widetilde {\Sigma^*}$ in $I^+(O)$ is thus a convex polyhedral surface with vertices $\left(p_\Delta\right)_{\Delta\in \mathfrak D}$  with 
 $\mathfrak D $ the set of connected component of $\widetilde \sing_0(M)$. 
 Then, the embedding of $\widetilde{\Sigma^*}$ is the boundary of the convex hull of the $(p_\Delta)_{\Delta\in \mathfrak D}$. 
 Finally, the Penner surface of $M$ is exactly the totally geodesic embedding of $\Sigma$ in $\E^{1,2}$ 
 induced by the suspension  of $\Sigma$ and thus  $$\pen(\susp_{\E^2}([\Sigma,S]))=[(\Sigma,S)].$$

\end{proof}

\subsection{Penner surface in non-linear Cauchy-compact spacetimes with BTZ}
The  aim of this section is to extend Theorem \ref{theo:teichlinear} to non-linear Cauchy-compact flat spacetimes with BTZ. We are not completely successful and 
only prove a partial result.
\begin{theo}\label{theo:teich_non_linear}
 Let $M$ be a Cauchy-compact Cauchy-maximal globally hyperbolic $\E^{1,2}_0$-manifold. 
 Let $(\Delta_i)_{i\in [1,s]}$ be the connected components of $\sing_0(M)$ 
 and let $(\bar p_i)_{i\in [1,s]}$ be a family of of points such that for all $i\in [1,s]$, $\bar p_i\in \Delta_i$.
 
 Then, there exists a unique convex polyhedral Cauchy-surface of $M$ with vertices $\bar p_1,\cdots,\bar p_s$.
 \end{theo}

Let $g \in \N$ and let $s>0$ such that $2g-2+s>0$. Let $[M,\bar p_1,\cdots, \bar p_s]$ be a point of  $\widetilde \M_{g,s}(\E^{1,2}_0)$,  the fundamental 
group of $M$ is a free group with at least 2 generators, it is thus non-abelian.
Therefore, $M$ automatically falls into case $(iv)$ of Theorem \ref{theo:barbot_mess}, in particular the linear part of the 
holonomy of $M$ is faithful and discrete. 
Let $G\subset \isom(\E^{1,2})$ be the image of the holonomy of $M$, let  $\widetilde \Omega$
be its augmented regular domain and let $G\subset \isom(\E^{1,2}) $ be the image of its holonomy. 
Let $(\Delta_i)_{i\in[1,s]}$ be a set of representative of the half lightlike line in $\widetilde \sing_0$, 
let $(p_i)_{i\in[1,s]}$ be the decoration of the $(\Delta_i)_{i\in[1,s]}$ and let $G_i=\stab(p_i)$. 

We still write 
$L:\isom(\E^{1,2})\rightarrow \SO_0(1,2)$ the projection on the linear part.

\begin{defi} Define $K(p)$, the closure of the convex hull of the set $\bigcup_{i=1}^s G p_i$ in $\widetilde \Omega$.
 
\end{defi}

\begin{lem} For all  $i\in [1,s]$,\label{lem:BTZ_K}
 $\Delta_i \cap K(p)= [p_i,+\infty[ $
\end{lem}
\begin{proof}
Let $i\in[1,s]$.
   \begin{itemize} 
  \item Assume $s\geq 2$, then take $p_j\neq p_i$. Since $p_j\in I^+(\Delta_i)$, there exists a unique $q\in \Delta_i$
  such that $ p_j\in \partial J^+(q)$. Moreover, $J^+(p_i)\cap \partial \Omega = ]p_i,+\infty[$ thus $q$ is in the past of $p_i$.
  Since  $G_i$ fixes point-wise $\Delta_i$, it acts linearly on $J^+(q)$.  From Corollary  \ref{cor:discrete}, $G p_j$ is discrete, 
  thus $G_i p_j$ is discrete. Let  $\phi\in G_i\setminus\{1\}$, we have $\phi_i^np_j= \lambda_n u_n$
  for some $\lambda_n\rightarrow +\infty$ and $u_n \rightarrow p_i$.
  For all $n\in \N$, the segment from $\phi^n p_j$ to  $p_i$ is a   subset of $K(p)$ and since $K(p)$ is closed, so is $[p_i,+\infty[$.
  
  Take a compact neighborhood $\U$ of $[q,p_j]$,  from Corollary  \ref{cor:discrete} the orbits $Gp_k$ for $k\in[1,s]$ are discrete, 
  thus $\bigcup_{k=1}^s Gp_k$ is discrete and there exists only finitely many points inside $\U$.  Moreover, $\bigcup_{k=1}^s G p_k \cap \Delta_i = p_i$, thus 
  $\U$ can be chosen small enough so that $\bigcup_{k=1}^s G p_k \cap \U = p_i$. Finally, $\Delta_i \cap K(p)= [p_i,+\infty[ $.
  \item Assume $s=1$ and $g\geq 1$. It is the same argument as before but instead of $p_j$ take $\psi p_i$ where $\psi$
  is some hyperbolic isometry of $G$.
 \end{itemize}
\end{proof}

\begin{lem}\label{lem:time_geod}
 Every timelike geodesic in $\widetilde \Omega$ intersects $\Sigma$. Moreover, once in $K(p)$, a future timelike geodesic 
 does not leave $K(p)$.
\end{lem}
\begin{proof}
Let $q\in \widetilde \Omega$ and let $u\in \H^2$. Applying some isometry of $\E^{1,2}$ we can assume $u=(1,0,0)$.
From Theorem \ref{theo:Teich}, the holonomy $\rho$ of $M$ is admissible.  The orbits 
any point $v\in\partial_\infty \H^2$ by $L(G)$ is thus dense. Let $u_1$ be a lightlike vector directing the BTZ line $\Delta$, then 
$L(G) u_1$ is dense and we can find $u_2$ and $u_3$ such that $u$ is in the interior of the convex hull of $u_1,u_2,u_3$ in 
$\H^2$. Let $q_1,q_2$ and $q_3$ be the respective decoration of the BTZ lines $\Delta_1,\Delta_2$ and $\Delta_3$ of direction $u_1,u_2$ and $u_3$ respectively.

For $T\in\R$, let $\Pi_T$ be the horizontal plane oh height $T$ and let $\pi$ be the vertical projection on $\Pi_0$.
From Lemma \ref{lem:BTZ_K}, for $T$ big enough $\Pi_T\cap K(p)$ contains the points $q_i+t_iu_i$, for $i\in\{1,2,3\}$ and
for some $t_i>T$. Furthermore, for $i\in\{1,2,3\}$ and $t>0$, $\pi(q_i+tu_i)=\pi(q_i)+t\pi(u_i)$ and since $u$ is in the interior 
of the convex hull of the $(u_j)_{j\in \{1,2,3\}}$ then the convex hull of  $\pi(\Delta_i\cap \Pi_T)$ is increasing with $T$
and their union for $T>0$ is $\Pi_0$. 
We deduce that for $t$ big enough, $q+t u \in K(p)$. 
This shows that $q+\R u$ intersects $K(p)$ and that for all $t$ such that $q+tu\in K(p)$, there exists $t'>t$ such that 
$q+t'u\in K(p)$. Since $K(p)$ is convex, the set of $t$ such that $q+tu\in K(p)$ is convex thus an interval and thus 
an interval of the 
form $[t_0,+\infty[$.

\end{proof}
\begin{cor} $K(p)$ is futures complete : $J^+(K(p))=K(p)$.
\end{cor}

\begin{proof} From Lemma \ref{lem:time_geod}, $I^+(K(p))=K(p)$ but $K(p)$ is closed and for all $q\in \E^{1,2}$, 
the closure of $I^+(q)$ is $J^+(q)$. Then $J^+(K(p))=K(p)$.
\end{proof}
\begin{cor}\label{cor:future_interior}
    For all $x \in \Omega \cap K(p), ~~J^+(x)\setminus \{x\} \subset Int(K(p))$,
\end{cor}
\begin{proof}
 Let $x\in \Omega\cap K(p)$,  using the same argument as for Lemma \ref{lem:time_geod}, we can find three points 
 $q_1,q_2,q_3$ on some BTZ-line such that the vertical projection of $x$ on the horizontal plane is in the convex 
 hull of the projections of $q_1,q_2$ and $q_3$.
 The facets $[q_iq_jx]$ are spacelike and  thus $J^+(x)$ is in the interior of the future of the convex hull of $\{q_1,q_2,q_3,x\}$.
 
\end{proof}

\begin{lem} $\partial K(p) $ is a closed achronal topological surface.
\end{lem}
\begin{proof}
 For $q\in K(p)$, let $f_q : \R^2 \rightarrow \R$ the function such that $\{(x,f(x)) : x\in \R^2 \} = \partial J^+(q)$. 
 For all $q\in K(p)$, $f_q$ is 1-lipschitz and, up to some global isometry, we can assume 
 $\forall q\in K(p), f_q\geq 0$,  thus $f:x\mapsto \inf_{q\in K(p)} f_q(x) $ is  well defined and 1-lipschitz.
 Since for all $q\in K(p)$,  $J^+(q) \subset K(p)$,  then $$\forall q\in K(p),\forall x\in \R^2,~~~~ (x,f(x)) \in K(p)$$
 Since $K(p)$ is closed then $\forall x\in \R^2, (x,f(x)) \in K(p)  $ and thus 
 $\partial_{\E^{1,2}} K(p) = \{(x,f(x)) : x \in \R^2\}$.  Finally, $\partial_{\E^{1,2}} K(p)$ is the graph
 of a 1-lipschitz function defined over $\R^2$ and thus $\partial_{\E^{1,2}} K(p)$ is a closed achronal topological surface.
\end{proof}

\end{proof}

\begin{prop}\label{prop:causal_curve}
 Every causal curves in $\Omega$ intersects $\partial K(p)$.
\end{prop}

\begin{proof}
 Define $\Omega'$ the interior Cauchy-developpement of $\partial K(p)$.
 Note that $\partial K(p)$ is not a lightlike plane, then from Theorem 1.5.1 p115 in \cite{algebre_2010}, 
 $\Omega'$ is a non-empty regular domain which is either future-complete, past-complete or between two lightlike planes.
  Since $\Omega'\subset \Omega$, it cannot be past-complete and from Lemma \ref{lem:time_geod} $\Omega'$ is future complete.

From Lemma \ref{lem:BTZ_K},  $\bigcup_{i=1}^s G]p_i,+\infty[\subset \partial K(p)$.
Therefore, the line of fixed point of the parabolic isometries of $G$ intersects the boundary of $\Omega'$ and thus,
by Theorem \ref{theo:max_BTZ_ext}, $\Omega'/G$ satisfies $(iii)$ of Lemma \ref{lem:abs_max_compact}. $\Omega'/G$ is thus 
absolutely maximal and $\Omega'=\Omega$.

\end{proof}

\begin{lem} \label{lem:boundary_K}
 The intersection of the boundary of $K(p)$ with the boundary of $\Omega$ is the union 
 of the future half BTZ-lines above $\phi p_i$ for $i\in [1,s]$ ans $\phi \in G$.
 $$\partial K(p) \cap \partial \Omega=\bigcup_{i=1}^s G[p_i,+\infty[ $$
\end{lem}
\begin{proof}
 Let $q \in \partial K(p) \cap \partial \Omega$ and  let $(q_n)_{n\in \N}$ be a sequence of points of the convex 
 hull $K(p)$ such that $q_n\xrightarrow{n\rightarrow +\infty} q$. 
 Write $$\forall n\in\N, \quad q_n = \sum_{\phi \in G}\sum_{i=1}^s \alpha_{i,\phi}^{(n)} \phi p_i$$
 
 There exists a past lightlike vector $u$ such that a plane of direction orthogonal to $u$ is a support plane of $\Omega$ 
 at $q$. Then, for all $x\in \Omega, \langle x |u\rangle > \langle q|u\rangle := r_0$ and for all $n\in \N$ :
    \begin{eqnarray*}
 \langle q_n |u\rangle&=&   \sum_{\phi \in G}\sum_{i=1}^s \alpha_{i,\phi}^{(n)} \langle \phi p_i|u\rangle \\
 &=&\sum_{(\phi,i) \in A_n}  \alpha_{i,\phi}^{(n)} \langle \phi p_i|u\rangle +\sum_{(\phi,i) \in B_n}  \alpha_{i,\phi}^{(n)} \langle \phi p_i|u\rangle
 \quad \quad \mathrm{With}~ \left\lbrace 
 \begin{array}{l}
A_n=\{(\phi,i) ~| ~\langle \phi p_i|u\rangle \leq r_0+\varepsilon_n  \} \\
      B_n=\{(\phi,i) ~|~ \langle \phi p_i|u\rangle > r_0+\varepsilon_n  \}    
      \\ \varepsilon_n \xrightarrow{n\rightarrow +\infty} 0
                     \end{array}\right.\\
 &\geq & \sum_{(\phi,i) \in A_n}  \alpha_{i,\phi}^{(n)}r_0 + \sum_{(\phi,i) \in B_n}  \alpha_{i,\phi}^{(n)} (r_0+\varepsilon_n)\\ 
 &\geq & r_0+ \varepsilon_n\sum_{(\phi,i) \in B_n}  \alpha_{i,\phi}^{(n)} 
 \end{eqnarray*}
Take $\varepsilon_n = \sqrt{\langle q_n |u\rangle -r_0}$ so that 
$$\sum_{(\phi,i) \in B_n}  \alpha_{i,\phi}^{(n)} \xrightarrow{n\rightarrow +\infty} 0.$$

Choose a past timelike vector $v$, there exists a support plane of $\Omega$ of direction $v^\perp$. Then 
for all $x\in \Omega, \langle x|v\rangle > r_1$ for some $r_1\in \R$.  For all $R>0$ and 
all $\varepsilon>0$, the domain
$$H_{R,\varepsilon}:=\left\{x\in \Omega ~\left|~\langle x|v\rangle <R \mathrm{~and~}\langle x|u\rangle < r_0+\varepsilon\right.\right\}$$
is relatively compact. Since $\bigcup_{i=1}^s Gp_i$ is discrete, for every $(R,\varepsilon)\in \R_+^2$ only finitely many points of  $\bigcup_{i=1}^s Gp_i$ are in $H_{R,\varepsilon}$.

If $u$ is not the direction of a BTZ-line,
then we can choose  $R$ such that $R>\langle q | v\rangle$ and $\varepsilon>0$ such that $H_{R,\varepsilon}$  contains
no decoration. If there exists a decoration point $a$ in $ (q+u^\perp)\setminus J^-(q)$, 
 $a-q$ is spacelike, then we choose $v$ such that
$\langle v | a-q\rangle > 0 $ and $R\in \R$ such that $\langle v | a\rangle>R>\langle v | q\rangle$. 
Again we can choose $\varepsilon$ small enough so that $H_{R,\varepsilon}$ contains no decoration point.
Either way, for $n\in\N$ such that $\varepsilon_n \leq \varepsilon$ : 
\begin{eqnarray}
\langle q_n|v \rangle &=&\sum_{(\phi,i) \in A_n}  \alpha_{i,\phi}^{(n)} \langle \phi p_i|v\rangle
+\sum_{(\phi,i) \in B_n}  \alpha_{i,\phi}^{(n)} \langle \phi p_i|v\rangle\\
&\geq & \sum_{(\phi,i) \in A_n}  \alpha_{i,\phi}^{(n)} R
+\sum_{(\phi,i) \in B_n}  \alpha_{i,\phi}^{(n)} r_1\\
\langle q|v \rangle+o(1)&\geq& R+o(1)\\
\end{eqnarray}
This is absurd since $R>\langle q|v \rangle$.

Therefore, there exist a decoration point in $J^-(q)$ and $q \in [p_0,+\infty[$ for some decoration point $p_0$.
 
\end{proof}

\begin{defi}
 Define $\widetilde \Sigma$  the boundary of $K(p)$ for the BTZ topology on $\widetilde \Omega$
 and define 
 $\Sigma := \widetilde \Sigma /G$.
\end{defi}
\begin{rem} We could have replaced $K(p)$ by the closure of  the convex hull of the decoration points for the BTZ topology instead of the usual topology. 
These two topologies coincide on $\Omega$ thus the closures are the same on $\Omega$. Both closure are then future complete
and thus for all decoration point $p$, the future set $I^+(p)$ is in both closures. Since the closure of $I^+(p)$  is $J^+(p)$ for both 
topologies, we end up having $[p,+\infty[$  in both closures. 
Finally, we see that the two closures are identical.

\end{rem}

\begin{lem} \label{lem:Sigma} The surface $\widetilde \Sigma$ is exactly
$$\widetilde \Sigma= \left(\Omega\cap \partial K(p)\right) \cup \{G p_1,\cdots,Gp_s\}.$$ 
Moreover, $$\Sigma = \partial_M(K(p)/G) $$
\end{lem}
\begin{proof} To begin with,the BTZ topology is thinner than the usual topology  thus 
$\widetilde \Sigma \subset \partial K(p)$. Since the BTZ topology coincides with the usual topology 
on $\Omega$ and since the decoration points are in $K(p)$, then $\left(\Omega\cap \partial K(p)\right) \cup \{G p_1,\cdots,Gp_s\} \subset \widetilde \Sigma$.
Now, let $p_i$ be a decoration point and let $q\in ]p_i,+\infty[$. 
Any open neighborhood of $q$ contains a set of the form $\Diamond_{q_1}^{q_2}$ with $q_1\in ]p_i,q[$. 
Since $\Diamond_{q_1}^{q_2}\subset J^+(p_i)\subset K(p)$, the point $q$ admits a neighborhood included into $K(p)$. 
Then $q$ is not in $\widetilde \Sigma$.

The projection $\widetilde \Omega \rightarrow M$ is open, the second point follows.
\end{proof}

\begin{prop} \label{prop:spacelike_facet}
$\Sigma$ is a convex spacelike polyhedral surface with finitely many 2-facets. Moreover the 0-facets are $\bar p_1,\cdots, \bar p_s$.
 \end{prop}
\begin{proof}
 From Lemma \ref{lem:Sigma}, $\Sigma$ intersects $\sing_0(M)$ exactly at $\{\bar p_1,\cdots,\bar p_2\}$ and is the boundary of a convex future set. $\Sigma$ is thus convex. Notice that Corollary 27 p415 of \cite{oneil} is valid for a $\E^{1,2}_0$-manifold and thus 
 $\Sigma$ is a topological surface in $M$.

 We first write a decomposition of $\partial K(p)$.
 Let $\Pi$ be some support plane of $K(p)$, from Lemma \ref{lem:time_geod} $\Pi$ cannot contains a timelike geodesic.
 Therefore, $\Pi$ is either lightlike or spacelike. From Corollary \ref{cor:future_interior}, 
 a lightlike plane cannot intersect $\Omega\cap K(p)$ and from Corollary \ref{prop:causal_curve}, it cannot intersect $\Omega$.
 Finally, from Lemma \ref{lem:boundary_K}, $\Pi\cap K(p)= [q,+\infty[$ for some decoration point $q$.
 
Assume $\Pi$ is a spacelike support plan at some point $q\in \Omega\cap K(p)$.
 then $K(p)\subset J^+(\Pi)$. If $\Pi$ contains less than two decoration points,
 then it can be moved slightly to obtain a plane $\Pi'$ such that $\{Gp_1,\cdots,Gp_s\}\subset J^+(\Pi')$ and $q\in I^-(\Pi')$,
 absurd.
  Therefore, $\Pi$ contains at least an open edge $e=]\phi_1 p_i, \phi_2 p_j[$ for some $\phi_1,\phi_2 \in G$ and $i,j \in [1,s]$ 
  and clearly, $q \in e$. Let $u\in \Pi^\perp$ a unit future timelike vector, 
  $u$ lies in a convex subset $\gamma$ of a geodesics of $\H^2$. 
  From  Corollary \ref{cor:future_interior}, we obtain that $\gamma$ is compact. We see that as long as $\Pi$
  doesn't contains three decoration points, $u$ is in the interior of $\gamma$.  
  Assume $\Pi$ contains at least 3 decoration points. 
  By discreteness of the set of decoration point, $\Pi$ contains a finite number of decoration points $\{q_1,\cdots,q_n\}$, $n\geq 3$.
  The convex hull $C$ of $(q_1,\cdots,q_n)$ is a convex polygon in $\Pi$. For all edge $f$ on the boundary of $C$, the set of 
  support plane at $f$ is paramtrized by a non trivial geodesic segment in $\H^2$. Therefore, $\Pi\cap K(p)=C$.  
  
  Let $\Pi$ lightlike plane support of $\Omega$ at some decoration point $q$. We can slightly rotate $\Pi$ to obtain a support
  spacelike plane at $q$. Then the set of future unit vector normal to a spacelike support plane at $q$ is a non-empty convex subset $H \subsetneq \H^2$. 
  A plane corresponding to a point in the boundary of $H$ contains at least 2 decoration points and thus $q$ 
  is a vertex of a spacelike 1-facet or a 2-facet. 
  
  Finally, a point of $\partial K(p)$ is either  on a spacelike facet of vertices in $\{Gp_1,\cdots, Gp_2\}$ 
  or on an infinite lightlike ray $[q,+\infty]$. Furthermore, every point of $\{Gp_1,\cdots, Gp_2\}$  is a vertex 
  of some 1-facet or 2-facet.
  Therefore, $\widetilde \Sigma$ is a spacelike polyhedral surface with 0-facets $Gp_1,\cdots, Gp_2$. 
  And then  $\Sigma$ is a spacelike polyhedral surface with 0-facets $\bar p_1,\cdots, \bar p_2$. 
  
  Assume there are infinitely many 2-facets, thus there is an accumulation point of 2-facets 
   in $\reg(\Sigma)$ and thus   in $\widetilde \Sigma \cap \Omega$.  Let $q$ be an accumulation point of 2-facets 
   in $\widetilde\Sigma \cap \Omega$, from the discussion above $q$ is either on a 1-facet or a 2-facet and in both cases,
   admits a neighborhood in $\Sigma$ with at most two 2-facets. This constradicts the fact that $q$ is an accumulation point.

\end{proof}

\begin{prop}\label{prop:cauchy_surface}
  $\Sigma $ is a Cauchy-surface of $M$.
\end{prop}
\begin{proof}

 Consider an inextensible future causal curve $\bar c:\R\rightarrow M$.  From Lemma 13 in \cite{brunswic_btz_ext}, it decomposes into a connected BTZ part 
 $\Delta$ and a connected regular part $\bar c^0$ such that $\Delta$ is  in the past of $\bar c^0$. 
 Write $c:\R\rightarrow \widetilde \Omega$ a lift of $\bar c$.
 
 \begin{itemize}
  \item If $\Delta\neq \emptyset$,  let $t_0\in \R$ such that $\bar c(t_0)=\max(\Delta)$. 
  Then $c(t_0) \in \phi[p_i,+\infty[$ for some $i\in [1,s]$ and $\phi \in G$.
 Then, $\phi p_i\in \Delta$ and $c(t_1)=\phi p_i \in \widetilde\Sigma$ for some $t_1\leq t_0 $. 

 For $t\in ]-\infty,t_0]$, $c(t)$ is on the BTZ line through $c(t_1)$ and thus $\forall t \in -\infty,t_0]\setminus\{t_1\}, c(t)\notin \widetilde\Sigma$. 
 Let $t>t_0$, then $c(t_1)\in \Omega\cap K(p)$, and then from Corollary \ref{cor:future_interior} $\forall t'>t, c(t')\notin \partial K(p)$. 
 Then $\forall t>t_0, c(t)\notin \widetilde\Sigma$.
 Finally, $\forall t\neq t_1, c(t)\notin \widetilde \Sigma$. 
 
  \item If $\Delta=\emptyset$ then $c\cap \partial K(p)\subset \Omega\cap \partial K(p)$ thus $c\cap \partial K(p)=c\cap \widetilde\Sigma$.
 From Proposition \ref{prop:causal_curve}, 
     $c\cap \partial K(p) \neq\emptyset$, then let $t_0:= c^{-1}(\min(c\cap \partial K(p)))$.
     From Proposition \ref{prop:causal_curve}, $\forall t>t_0, c(t)\notin \partial K(p)$ and thus 
     $c \cap \partial K(p)=\{c(t_0)\}$.
 \end{itemize}
Therefore, $\Sigma$ is a Cauchy-surface of $M$. 
\end{proof}

\addtocounter{theo}{-1}

\begin{theo}
 Let $M$ be a Cauchy-compact Cauchy-maximal globally hyperbolic $\E^{1,2}_0$-manifold. 
 Let $(\Delta_i)_{i\in [1,s]}$ be the connected components of $\sing_0(M)$ 
 and let $(\bar p_i)_{i\in [1,s]}$ be a family of of points such that for all $i\in [1,s]$, $\bar p_i\in \Delta_i$.
 
 Then, there exists a unique convex polyhedral Cauchy-surface of $M$ with vertices $\bar p_1,\cdots,\bar p_s$.
 \end{theo}

 \begin{proof}
  Consider $\Sigma = \partial_M(\partial K(p)/G)$. By Proposition \ref{prop:cauchy_surface}, $\Sigma$ is a Cauchy-surface and by Proposition 
  \ref{prop:spacelike_facet}, $\Sigma$ is a convex polyhedral surface intersecting $\Delta_i$ exactly at $\bar p_i$.
  
  Let $\Sigma_1$ be another convex polyhedral Cauchy-surface of vertices $\bar p_1,\cdots, \bar p_s$. 
  On the one hand, $J^+(\Sigma_1)$ contains $K(p)$ and thus $\Sigma$ is in the future of $\Sigma_1$. 
  On the other hand, consider an edge of $\Sigma_1$. Since the vertices of $\Sigma_1$ are $\bar p_1,\cdots, \bar p_s$, 
  this edge is a geodesic segment from some $\bar p_i$ to some $\bar p_j$ and thus belongs to $K(p)/G$. 
  Finally, $\Sigma_1$ lies in the future of $\Sigma$ and these two Cauchy-surfaces are thus equal.
  
 \end{proof}

    \addcontentsline{toc}{section}{References}
    \bibliographystyle{alpha}
    \bibliography{note} 

\begin{thebibliography}{AGH98}

\bibitem[AGH98]{andersson_time_function}
Lars Andersson, Gregory~J. Galloway, and Ralph Howard.
\newblock The cosmological time function.
\newblock {\em Classical Quantum Gravity}, 15(2):309--322, 1998.

\bibitem[alg10]{algebre_2010}
{\em Algèbre, dynamique et analyse pour la géométrie : aspects récents}.
\newblock Ellipses, 2010.

\bibitem[Bar05]{barbot_globally_2004}
Thierry Barbot.
\newblock Globally hyperbolic flat spacetimes.
\newblock {\em {Journal of Geometry and Physics}}, 53, no.2:123--165, 2005.

\bibitem[BB09]{MR2499272}
Riccardo Benedetti and Francesco Bonsante.
\newblock Canonical {W}ick rotations in 3-dimensional gravity.
\newblock {\em Mem. Amer. Math. Soc.}, 198(926):viii+164, 2009.

\bibitem[BBS11]{Particles_1}
Thierry Barbot, Francesco Bonsante, and Jean-Marc Schlenker.
\newblock Collisions of particles in locally {A}d{S} spacetimes {I}. {L}ocal
  description and global examples.
\newblock {\em Comm. Math. Phys.}, 308(1):147--200, 2011.

\bibitem[BBS14]{Particles_2}
Thierry Barbot, Francesco Bonsante, and Jean-Marc Schlenker.
\newblock Collisions of particles in locally {A}d{S} spacetimes {II}. {M}oduli
  of globally hyperbolic spaces.
\newblock {\em Comm. Math. Phys.}, 327(3):691--735, 2014.

\bibitem[BC]{spinexp}
Thierry Barbot and Meusburger Catherine.
\newblock Particles with spin in flat spacetimes in expansion.
\newblock in preparation.

\bibitem[Bon03]{bonsante_flat_2003}
Francesco Bonsante.
\newblock Flat {Spacetimes} with {Compact} {Hyperbolic} {Cauchy} {Surfaces}.
\newblock November 2003.

\bibitem[Bru16]{brunswic_btz_ext}
L{\'e}o Brunswic.
\newblock {BTZ extensions of globally hyperbolic singular flat spacetimes}, May
  2016.
\newblock preprint : https://arxiv.org/abs/1605.05530.

\bibitem[CBG69]{MR0250640}
Yvonne Choquet-Bruhat and Robert Geroch.
\newblock Global aspects of the {C}auchy problem in general relativity.
\newblock {\em Comm. Math. Phys.}, 14:329--335, 1969.

\bibitem[CI01]{Troyanov}
M.~Troyanov H.~Clémençon C.~Indermitte, Th. M~Liebling.
\newblock Voronoi diagrams on piecewise at surfaces and an application to
  biological growth.
\newblock {\em Theoretical Computer Science}, (263):263--274, 2001.

\bibitem[EP88]{MR918457}
D.~B.~A. Epstein and R.~C. Penner.
\newblock Euclidean decompositions of noncompact hyperbolic manifolds.
\newblock {\em J. Differential Geom.}, 27(1):67--80, 1988.

\bibitem[Ger70]{MR0270697}
Robert Geroch.
\newblock Domain of dependence.
\newblock {\em J. Mathematical Phys.}, 11:437--449, 1970.

\bibitem[Gol]{goldman_course_note}
William~M. Goldman.
\newblock Projective geometry on manifolds.

\bibitem[Gol84]{MR762512}
William~M. Goldman.
\newblock The symplectic nature of fundamental groups of surfaces.
\newblock {\em Adv. in Math.}, 54(2):200--225, 1984.

\bibitem[Hit92]{MR1174252}
N.~J. Hitchin.
\newblock Lie groups and {T}eichm\"uller space.
\newblock {\em Topology}, 31(3):449--473, 1992.

\bibitem[Mes07]{mess}
Geoffrey Mess.
\newblock Lorentz spacetimes of constant curvature.
\newblock {\em Geom. Dedicata}, 126:3--45, 2007.

\bibitem[Min61]{MR0132379}
George~J. Minty.
\newblock On the maximal domain of a ``monotone'' function.
\newblock {\em Michigan Math. J.}, 8:135--137, 1961.

\bibitem[Min64]{MR0167859}
George~J. Minty.
\newblock On the monotonicity of the gradient of a convex function.
\newblock {\em Pacific J. Math.}, 14:243--247, 1964.

\bibitem[MS08]{MR2436235}
Ettore Minguzzi and Miguel S{\'a}nchez.
\newblock The causal hierarchy of spacetimes.
\newblock In {\em Recent developments in pseudo-{R}iemannian geometry}, ESI
  Lect. Math. Phys., pages 299--358. Eur. Math. Soc., Z\"urich, 2008.

\bibitem[O'N83]{oneil}
Barett O'Neill.
\newblock {\em Semi-Riemannian geometry}.
\newblock Academic Press, 1983.

\bibitem[Pen87]{penner1987}
R.~C. Penner.
\newblock The decorated teichmüller space of punctured surfaces.
\newblock {\em Comm. Math. Phys.}, 113(2):299--339, 1987.

\bibitem[Pen12]{MR3052157}
Robert~C. Penner.
\newblock {\em Decorated {T}eichm\"uller theory}.
\newblock QGM Master Class Series. European Mathematical Society (EMS),
  Z\"urich, 2012.
\newblock With a foreword by Yuri I. Manin.

\bibitem[Rat]{ratcliff_foundation}
John Ratcliffe.
\newblock {\em Foundations of Hyperbolic Manifolds}.
\newblock 149. Springer-Verlag New York.

\end{thebibliography}
\end{document}